\documentclass{conm-p-l-version} 
\usepackage{amsfonts}

\usepackage{xr}
\newcommand{\bfc}{{\Bbb C}}
\newcommand{\bfz}{{\Bbb Z}} 
\newcommand{\bfq}{{\Bbb Q}} 
\newcommand{\noi}{\noindent}
\newcommand{\tw}{T(E_w)}
\newcommand{\tv}{T(E_v)}
\newcommand{\G}{\Gamma}
\newcommand{\cvg}{{\mathcal V}(\Gamma)}
\newcommand{\evg}{{\mathcal E}(\Gamma)}
\newcommand{\cale}{{\mathcal E}}

\newcommand{\calc}{{\mathcal C}}
\newcommand{\calg}{{\mathcal G}}
\newcommand{\calw}{{\mathcal W}}
\newcommand{\calv}{{\mathcal V}}
\newcommand{\cala}{{\mathcal A}}

\newcommand{\caly}{{\mathcal Y}}
\newcommand{\bn}{{\mathbf n}}
\newcommand{\bd}{{\mathbf d}}
\newcommand{\fedog}{{\mathcal G}(\Gamma, (\bn,\bd))}
\newcommand{\block}{B(n_v,n_w,d_e)}
\newcommand{\bgu}{B(\Gamma,\bn)}
\newcommand{\agn}{A(\Gamma,\bn)}
\newcommand{\coker}{\mbox{coker}\,}
\newcommand{\finv}{f^{-1}(0)}
\newcommand{\xfn}{X_{f,N}\,}
\newcommand{\gxf}{G(X,f)}
\newcommand{\Gxf}{\Gamma (X,f)}
\newcommand{\fegyz}{\{ f(x,y)=z^{N}\} }
\newcommand{\geo}{h_{geom}}
\newcommand{\szam}[2]{\#_{#1}^{#2}(f)}

\begin{document}

\title[Resolution Graphs of Cyclic Coverings]{Resolution Graphs of Some Surface Singularities, I.\\
(Cyclic Coverings)}
\author{ Andr\'as N\'emethi}
\address{Department of Mathematics, The Ohio State University, Columbus,
Ohio 43210}
\email{nemethi@math.ohio-state.edu}

\begin{abstract}
The article  starts with some introductory material about
resolution graphs of normal surface singularities (definitions,
topological/homological properties, etc). Then we  discuss the problem of
$N$-cyclic coverings $(X_{f,N},0)\to (X,0)$ of $(X,0)$,
branched along $(\{f=0\},0)$, where
$f:(X,0)\to ({\bfc},0)$ is the germ of an analytic function.
We present non--trivial examples in order to show that from the embedded
resolution graph $\Gamma(X,f)$ of  $f$  it is not possible to recover the
resolution graph of $(X_{f,N},0)$. The main results are the construction of
a ``universal covering graph'' of $\Gamma(X,f)$ from the topology of the
germ $f$, and the completely combinatorial construction  of the resolution
graph of $(X_{f,N},0)$ from this universal graph of $f$ and the integer $N$.
For this we also prove some classification theorems of ``graph coverings'',
results which are purely graph-theoretical. In the last part, we connect the
properties of the universal covering graph with the topological invariants
of $f$, e.g. with the nilpotent part of its algebraic monodromy.
\end{abstract}

\maketitle

\subsection*{Introduction.} 

The present article has several goals. First of all, it is an article
with some  expository character, presenting several aspects of the resolution
graphs of normal surface singularities. Moreover, it creates the good language 
and right point of view in the graph--codifications of important geometric
constructions, such as cyclic coverings of surface singularities, series of
singularities, degeneration of curves, etc. This is realized by the 
introduction of a new combinatorial object: the 
 ``covering of graphs'' (Section \ref{sec1}). 
We  will present 
two main applications: the case of cyclic coverings of normal surface 
singularities (the present article, in the  sequel: Part I), 
and the case of some series of singularities 
(a joint work with \'Agnes Szil\'ard \cite{partII}, in the sequel: Part II). 

Consider a normal surface singularity $(X,0)$ and   a germ of an analytic 
function  $f:(X,0)\to (\bfc,0)$. Let $X_{f,N}$ be the (normalized)
cyclic $N$--covering of $(X,0)$ branched along $(\{f=0\},0)$. 
The final goal of Part I is to recover the resolution graph 
$\G(X_{f,N})$
of $X_{f,N}$  from some invariants of the pair $(X,f)$ and the integer $N$. 
Particular cases suggest that  $\G(X_{f,N})$
might be computable from the embedded resolution graph $\G(X,f)$ of $f$
and the integer $N$, but this is not true in general. In fact, $\G(X,f)$
 and $N$ codify all the local data of the covering, but they miss some global
information. This global data is in a close relationship with the 
monodromy representation $arg_*$ of the Milnor fibration associated with $f$,
which plays an 
important role in the global geometry  of all the cyclic coverings. 

In the second section we present the basic definitions and some of the 
(elementary) properties of  embedded resolution graphs.  This section also
serves  as an introduction to  Part II. Then we discuss ``local/global''
properties, showing  e.g. that the embedded resolution graph does not 
determine the  representation $arg_*$.
The additional global information of $arg_*$ will be 
codified  in another graph $G(X,f)$, called the ``universal covering graph''
of $\G(X,f)$. 

This  already motivates the material of Section \ref{sec1}, 
where  a 
detailed presentation of the ``cyclic covering of graphs'' can be found. 
The section 
culminates in the classification theorems, which are indispensable  in  the 
applications. 

In the third section we present an algorithm which provides
 the resolution graph $\G(X_{f,N})$,
in a purely combinatorial way, from the universal covering  graph $G(X,f)$
of $\G(X,f)$ and the integer $N$. In fact, each graph
$\G(X_{f,N})$ appears as a graph covering 
$\G(X_{f,N})\to \G(X,f)$ modified by  Hirzebruch--Jung strings. 
The family of graphs $\{\G(X_{f,N})\}_N$ 
behaves  like a  complicated series of
graphs, which is coordinated by a unique graph, the universal
covering of $\G(X,f)$. The global nature of the construction is
emphasized by non--trivial examples. On the other hand, we always stress
the particular cases when the global information can be recovered from the 
local (i.e. when $G(X,f)$ is determined from $\G(X,f)$), 
e.g. when the link of $(X,0)$ is a rational homology sphere. 

Each section and subsection has its own introduction, where the reader can 
find more guiding information (see also the introduction of Part II).

\section{${\bfz}$--coverings of graphs}\label{sec1}

\bekezdes{}{a} In this section  we present some graph--theoretical 
constructions.  The construction
 of ``cyclic covering of graphs'' is motivated by 
the structure of the resolution
graphs of cyclic coverings of the normal surface singularities
$(X,0)$ with branch locus $(f^{-1}(0),0)\subset (X,0)$
where $f:(X,0)\to (\bfc,0)$ is a germ of an analytic function
(cf. \ref{1.10}).
All these resolution graphs are controlled by a single ``universal
covering graph''
which is a ${\bfz}$--covering of the embedded resolution graph 
of the pair $(f^{-1}(0),0)\subset (X,0)$.

Graph coverings can be defined for an arbitrary discrete group $G$,
and they have important applications even for $G\not={\bfz}$.
The general description, together with some applications, will be 
published elsewhere. In this article we only discuss the case 
$G={\bfz}$, and we call the corresponding coverings ``cyclic'' or
``${\bfz}$--coverings''.

\bekezdes{Notations.}{b} For any graph $\Gamma$, we denote the 
set of vertices  by 
$\cvg$ and the set of edges by $\evg\subset\cvg\times\cvg$.
If there is no danger of confusion, we denote them simply by
${\calv}$ and ${\cale}$.

For simplicity, we will assume that our graphs have no loops.
The interested reader can
easily extend all our results  for graphs with loops.
If $v\in{\calv}$ is a fixed  vertex, let ${\cale}_v$ be the set of edges 
which have  $v$ as one of their endpoints. 
If $|\Gamma|$ is the topological realization of the graph $\Gamma$, then we 
denote  the rank of $H_1(|\Gamma|,\bfz)$ by $c_{\Gamma}$, i.e. $c_{\Gamma}$
is the number of independent cycles of the graph $\Gamma$.

\bekezdes{Definitions.}{4.1}

{\em A morphism of graphs $p:\Gamma_1\rightarrow\Gamma_2$ consists of two maps
$p_{\calv}: {\calv}(\Gamma_1)\rightarrow {\calv}(\Gamma_2)$ and
$p_{\cale}: {\cale}(\Gamma_1)\rightarrow {\cale}(\Gamma_2)$, such that
if $e\in{\cale}(\Gamma_1)$ has endpoints $v_1$ and $v_2$,
then $p_{\cale}(e)$ is an edge
in $\Gamma_2$ with endpoints  
$p_{\calv}(v_1)$ and  $p_{\calv}(v_2)$.
If $p_{\calv}$ and $p_{\cale}$ are isomorphisms of sets, then we say that 
$p$ is an isomorphism of graphs.

If $\Gamma$ is a graph, we say that ${\bfz}$ acts on $\Gamma$, if there
are group--actions $a_{\calv}: {\bfz}\times\calv\rightarrow\calv$ and 
 $a_{\cale}: {\bfz}\times\cale\rightarrow\cale$ of \ ${\bfz}$ 
with the following  compatibility property: 
if $e\in {\cale}$ has endpoints $v_1$ and $v_2$, then
$a_\cale (h,e)$ has endpoints $a_\calv (h,v_1)$ and $a_\calv (h,v_2)$. 
The action is trivial if $a_\calv$ and $a_\cale$
are trivial actions. 

If \ ${\bfz}$ acts on both $\Gamma_1$ and $\Gamma_2$, then  a morphism 
$p:\Gamma_1\to \Gamma_2$ is equivariant if the maps $p_{\calv}$ and
$p_{\cale}$ are equivariant with respect to the actions of ${\bfz}$. If 
additionally $p$ is 
an isomorphism then it is called an equivariant isomorphism of graphs.}

\bekezdes{-- The ``segment graph''.}{c}
The simplest graph is the 
 ``segment graph'' $S$ which has two vertices 
$\calv=\{ v_1, v_2\}$ and one edge 
$\cale =\{e=(v_1, v_2)\}$. In the sequel  $S$ will always carry the
 trivial action of the structure group ${\bfz}$.

\bekezdes{-- The ``standard blocks'' of $\bfz$--coverings.}{4.2}\ 
A ``standard block'' $B$, which covers the segment graph  $S$,
 can be constructed as follows.

Fix three strictly positive integers $n_1, n_2$ and $d$, and set
 $[n,m]=l.c.m.\{n,m\}$.
The standard block  $B(n_1, n_2, d)$ is a graph which
consists of $n_1+n_2$ 
vertices $\calv = \{ P_1^1,...,P_{n_1}^1, P_1^2,...,P_{n_2}^2\}$ 
and $d\cdot [n_1,n_2]$ edges $\cale = \{ e_1,...,e_{d[n_1,n_2]}\}$
and has the following structure:
for each $k\in \{ 1,\ldots, d[n_1,n_2]\}$, 
the endpoints of the edge $e_k$ are $P_{i(k)}^1$ and  $P_{j(k)}^2$, where $i(k)\equiv k 
\ (mod\ {n_1})$ and  $j(k)\equiv k \ (mod\ {n_2})$. Notice that  $P_{i}^1$ and  $P_{j}^2$ are connected
by exactly $d$ edges. The graph  $B(n_1, n_2, d)$ has a natural $\bfz$--action
given by $a_\cale (n,e_k)=e_l$, where $l\equiv k+n \ (mod\ {d[n_1,n_2]})$ and
$a_\cala (n,P_k^i)=P_l^i$, where $l\equiv k+n \ (mod\ {n_i})\ \ (i=1,2)$.

Consider the ``segment graph'' $S$ and the trivial $\bfz$--action on it.
Then $p: B(n_1,n_2,d)\rightarrow S$ defined by 
$p(P_j^i)=v_i\ \ (i=1,2)$ and $p(e_k)=(v_1,v_2)$ is an equivariant morphism.

Actually, almost all the finite
coverings of $S$ are ``standard blocks''. Indeed,
consider an arbitrary  equivariant morphism $p:B\to S$ of graphs,
where $B$ is a finite connected graph with ${\cale}\not=\emptyset$ and
with a $\bfz$--action, 
$p$ is an equivariant morphism such that
the restriction of the action of $\bfz$ on 
$p^{-1}(v_i)$ ($i=1,2$), 
respectively on 
$p^{-1}((v_1,v_2))$, is transitive. Then, it is not difficult to show
that $p:B\to S$ is equivalent with the standard block $B(n_1,n_2,d)$,
where $n_i\bfz$ (respectively $d[n_1,n_2]\bfz$) is the maximal 
subgroup of $\bfz$
which acts trivially on $p^{-1}(v_i)$ ($i=1,2$) (respectively on 
 $p^{-1}((v_1,v_2))$). 

\bekezdes{-- Cyclic (or $\bfz-$) coverings of graphs.}{4.3}\
Now, consider an arbitrary graph $\Gamma$. 
We assume that $\bfz$ acts on $\Gamma$ in a trivial way. 
Any edge $e$ with endpoints $v_1$ and $v_2$  determines a natural 
``segment subgraph'' $S_e = (\{v_1,v_2\},\{e\} )$ of $\Gamma$.

\bekezdes{Definition.}{def1}\ {\em A $\bfz$--covering, or cyclic covering,
 of the finite graph $\Gamma$ consists of a 
finite graph $G$, that carries a $\bfz$--action, 
together with an equivariant morphism 
$p:G\rightarrow\Gamma$ such that the restriction of the $\bfz$--action on
any set of type $p^{-1}(v)$ ($v\in \cvg$), respectively $p^{-1}(e)$ ($e\in\evg$), is transitive.}

\vspace{2mm}

This definition can be reformulated in terms of ``standard blocks'' as follows.
For any $v\in\cvg$ (respectively
 edge $e\in\evg$ with endpoints $\{v_1,v_2\}$), let 
$n_v\bfz$ (respectively $d_e[n_{v_1},n_{v_2}]\bfz$) be  the maximal 
subgroup of $\bfz$
which acts trivially on $p^{-1}(v_i)$ ($i=1,2$) (respectively on 
 $p^{-1}(e)$). This defines a 
 system of strictly positive integers
$(\bn, \bd)=\left\{ \{n_v\}_{v\in\calv(\Gamma)}; \{d_e\}_{e\in\cale(\Gamma)}
\right\}$. This system will be called {\em covering data}. 

\bekezdes{Definition.}{4.4}\ {\em 
Fix a graph $\Gamma$ with a trivial $\bfz$--action, and a system of integers
(covering data) 
$(\bn, \bd)=\left\{ \{n_v\}_{v\in\calv}; \{d_e\}_{e\in\cale}\right\}$.
A $\bfz$--covering (or cyclic covering) of type $(\bn,\bd)$ 
 of the graph $\Gamma$ consists of a 
graph $G$, that carries a $\bfz$--action, together with an equivariant morphism 
$p:G\rightarrow\Gamma$ such that for any $v\in\cvg$ (resp. $e\in\evg$) the set
$p^{-1}(v)$ (resp. $p^{-1}(e)$) consists of $n_v$ vertices (resp. 
$n_e=d_{e}\cdot [n_{v_1},n_{v_2}]$ edges).
Moreover, we assume that for any edge $e$, the subgraph
$p^{-1}(S_e)$ is a ``standard block'' 
$B=B(n_{v_1},n_{v_2}, d_e)$ such that the restriction of the $\bfz$--action
of $\bfz$ to $p^{-1}(S_e)$ coincides with the natural $\bfz$--action of $B$.}

\vspace{2mm}

The above definition shows that 
any cyclic covering $G\rightarrow\Gamma$ is constructed from 
standard blocks which are  glued equivariantly along the vertices 
$\{p^{-1}(v)\}$, where
$v$ runs over the vertices of $\Gamma$ with degree $\geq 2$.

\bekezdes{Definition.}{def2} {\em 
Two cyclic coverings $p_i : G_i\rightarrow\Gamma\ \ (i=1,2)$ are equivalent $(G_1 \sim G_2)$
if there is an equivariant isomorphism $q:G_1\rightarrow G_2$ such that $p_2 \circ q=p_1$.

The set of equivalence classes of cyclic coverings of $\Gamma$, associated with a  system of 
integers $(\bn,\bd)$, is denoted by $\fedog$.}

\bekezdes{Examples.}{4.5}\ 

a.) If $\Gamma=S$, then $\fedog$ has exactly one element for any
$(\bn,\bd)$.

b.) Let $\Gamma$ be the cyclic graph with two vertices and two edges:
$\cvg =\{ v_1,v_2\}$, $\evg = \{ e_1,e_2\}$ with both $e_i\ \ (i=1,2)$ having
$v_1$ and $v_2$ as endpoints.

Set $n_1=n_2=n$. Then $\fedog$
has exactly $n$ elements. Notice that we can have graphs that are not
 equivalent as cyclic coverings over $\Gamma$, but they are isomorphic as 
graphs. (Take e.g. the case $n=3$.)

The fact that in this case  $\ \fedog$ is independent of the choice of $\bd$, 
is not a particularity of this example: it is true in the most general situation, 
cf. (\ref{4.10}).

\bekezdes{-- The trivial covering of  $(\Gamma, (\bn,\bd))$.}{4.6}

There is a special element in $\fedog$ which can be 
constructed as follows. Fix a distinguished edge  in each standard block 
$B(n_{v_1},n_{v_2}, d_e)$. Then construct $G$ 
in such a way that whenever the vertex $v\in\cvg$ of $\Gamma$ is adjacent to 
the edges $\{e\}_{e\in\cale_v}$ of $\Gamma$, 
then the distinguished edges of all the blocks
$\{ p^{-1}(e)\}_{e\in\cale_v}$ have a common endpoint 
(which is one of the vertices
of $G$ in $p^{-1}(v)$).
Notice that this condition together with the existence of the $\bfz$--action
determines all the other adjacency relations.

The equivalence class of $G$ constructed in this way does not depend on the choice
of the distinguished edges in the standard blocks. Indeed, if we have two 
different choices of the distinguished edges in the standard blocks which
provide two graphs
$G_1$ and $G_2$ by the above construction  (both covering $\Gamma$), 
then for any edge $e$ of
 $\Gamma$ denote  the distinguished edges in the blocks $p_i^{-1}(S_e)$
by $e_i$ \ \ $(i=1,2)$.
Then the map $q(e_1)=e_2$ (for any $e$) can be extended to an
equivariant isomorphism $q: G_1\rightarrow G_2$ with $p_2 \circ q=p_1$.

The covering just constructed (or its class) is denoted by $p:T\rightarrow\Gamma$ and it
is called the {\em trivial cyclic covering of $(\Gamma, (\bn,\bd))$. }

The trivial covering is  characterized 
 by the existence of a (non-equivariant!) morphism
of graphs  $s:\Gamma\rightarrow T$ with
$p\circ s=id_\Gamma$ (i.e. $s$ is a section with $s(e)=$distinguished edge above $e$).

\bekezdes{-- $\fedog$ as a homogeneous space. Classification.}{4.7}\

If $p:G\rightarrow\Gamma$ is a cyclic covering of $\Gamma$, 
then $G$ is obtained by
an equivariant gluing of $\#\evg$ standard blocks. Above any vertex
$v\in\cvg$ we have to 
glue together $\#{\cale}_v$ standard blocks. 

Regard $\Gamma$ as a union of 
segments, with some  endpoints glued together. Then it is useful to 
introduce an index set of all the endpoints of the segments. This is:
$I= \coprod\limits_{v}{\cale}_v = \coprod\limits_{v} \{e_v\}_{e_v\in
{\cale}_v}$.

Now consider the following group indexed exactly over this set:
$B(\Gamma, \bn ) = 
\prod\limits_{v}\prod\limits_{e_v\in{\cale}_v} {\bfz}_{n_v}=
\prod\limits_{e_v\in I}{\bfz}_{n_v}$.
Then 
$B(\Gamma, \bn )$ describes the equivariant gluings in $G$. A typical element in
$B(\Gamma, \bn )$ is $\{ b_{e_v}\}_{e_v\in I}$, where $b_{e_v}\in 
{\bfz}_{n_v}$.
A generator set $\{  g_{e_v}\}_{e_v\in I}$ has the form 
$ g_{e_v}=(0,...,\hat{1},...0)$, where  all the entries are zero except 
the place $e_v$, where we put the generator $\hat{1}$ of ${\bfz}_{n_v}$.\\

$B(\Gamma, \bn )$ acts in a natural way on the set $\fedog$. We give the action of 
$ g_{e_v}$ for each $e_v\in I$. For this, fix $v\in\cvg$ and an edge $e=e_v$
with endpoints $v$ and $w$. 
If $G\rightarrow\Gamma$ is a covering, then above the segment $S_e$ we have the block
$B(n_v,n_w,d_e)$. In particular, above the vertex $v$ we have $n_v$ vertices cyclically 
permuted by the $\bfz$--action. Call them $P_1,...,P_{n_v}$ so that the action is
$P_1\rightarrow P_2\rightarrow ...\rightarrow P_{n_v}\rightarrow P_1$. We can detach the 
endpoint $v$ of the segment $e$ from the graph $\Gamma$.
In the graph $G$ this means that we separate the endpoints of the block $\block$, that
stay above $v$,  from $G$. In this way we obtain the graph $\tilde{G}_v$.
It is represented in the right hand side of the next diagram.

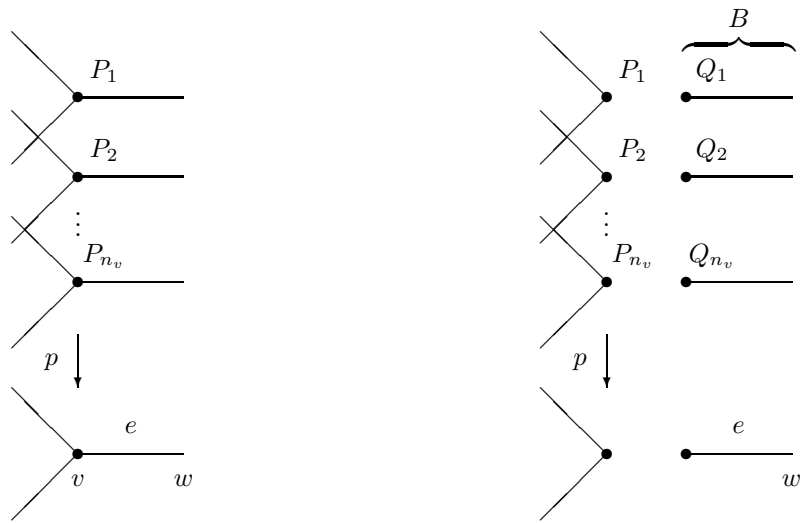
\begin{picture}(400,200)(0,0)

\put(50,25){\circle*{4}}
\put(50,90){\circle*{4}}
\put(50,130){\circle*{4}}
\put(50,160){\circle*{4}}
\put(50,25){\line(1,0){40}}
\put(50,90){\line(1,0){40}}
\put(50,130){\line(1,0){40}}
\put(50,160){\line(1,0){40}}
\put(50,25){\line(-1,-1){25}}
\put(50,90){\line(-1,-1){25}}
\put(50,130){\line(-1,-1){25}}
\put(50,160){\line(-1,-1){25}}
\put(50,25){\line(-1,1){25}}
\put(50,90){\line(-1,1){25}}
\put(50,130){\line(-1,1){25}}
\put(50,160){\line(-1,1){25}}
\put(50,70){\vector(0,-1){20}}
\put(40,60){\makebox(0,0){$p$}}
\put(60,100){\makebox(0,0){$P_{n_v}$}}
\put(60,140){\makebox(0,0){$P_2$}}
\put(60,170){\makebox(0,0){$P_1$}}
\put(50,115){\makebox(0,0){$\vdots$}}
\put(50,15){\makebox(0,0){$v$}}
\put(90,15){\makebox(0,0){$w$}}
\put(70,35){\makebox(0,0){$e$}}

\put(250,25){\circle*{4}}
\put(250,90){\circle*{4}}
\put(250,130){\circle*{4}}
\put(250,160){\circle*{4}}
\put(280,25){\circle*{4}}
\put(280,90){\circle*{4}}
\put(280,130){\circle*{4}}
\put(280,160){\circle*{4}}
\put(280,25){\line(1,0){40}}
\put(280,90){\line(1,0){40}}
\put(280,130){\line(1,0){40}}
\put(280,160){\line(1,0){40}}
\put(250,25){\line(-1,-1){25}}
\put(250,90){\line(-1,-1){25}}
\put(250,130){\line(-1,-1){25}}
\put(250,160){\line(-1,-1){25}}
\put(250,25){\line(-1,1){25}}
\put(250,90){\line(-1,1){25}}
\put(250,130){\line(-1,1){25}}
\put(250,160){\line(-1,1){25}}
\put(250,70){\vector(0,-1){20}}
\put(240,60){\makebox(0,0){$p$}}
\put(260,100){\makebox(0,0){$P_{n_v}$}}
\put(260,140){\makebox(0,0){$P_2$}}
\put(260,170){\makebox(0,0){$P_1$}}
\put(290,100){\makebox(0,0){$Q_{n_v}$}}
\put(290,140){\makebox(0,0){$Q_2$}}
\put(290,170){\makebox(0,0){$Q_1$}}
\put(250,115){\makebox(0,0){$\vdots$}}
\put(320,15){\makebox(0,0){$w$}}
\put(300,35){\makebox(0,0){$e$}}
\put(300,180){\makebox(0,0){$\overbrace{\hspace{1.5cm}}$}}
\put(300,190){\makebox(0,0){$B$}}
\end{picture}

Hence, we can re-obtain $G$ from $\tilde{G}_v$
if we re-glue $P_i$ and $Q_i$ $(i=1,...,n_v)$.
We emphasize that the action is $P_i\rightarrow P_{i+1}$ and similarly 
$Q_i\rightarrow Q_{i+1}$.
Now, by definition,
$g_{e_v} * G$ is obtained by the equivariant gluing $P_1=Q_2,\ P_2=Q_3,\ldots,
P_{n_v}=Q_1$.

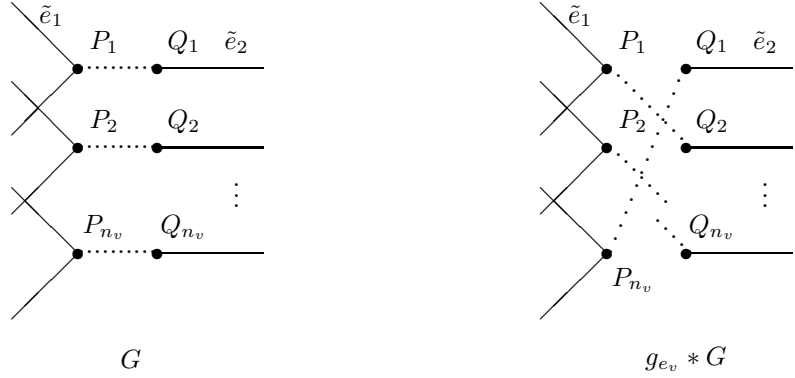
\begin{picture}(400,160)(0,40)
\put(50,90){\circle*{4}}
\put(50,130){\circle*{4}}
\put(50,160){\circle*{4}}
\put(80,90){\circle*{4}}
\put(80,130){\circle*{4}}
\put(80,160){\circle*{4}}
\put(110,115){\makebox(0,0){$\vdots$}}
\put(80,90){\line(1,0){40}}
\put(80,130){\line(1,0){40}}
\put(80,160){\line(1,0){40}}
\put(50,90){\line(-1,-1){25}}
\put(50,130){\line(-1,-1){25}}
\put(50,160){\line(-1,-1){25}}
\put(50,90){\line(-1,1){25}}
\put(50,130){\line(-1,1){25}}
\put(50,160){\line(-1,1){25}}
\put(60,100){\makebox(0,0){$P_{n_v}$}}
\put(60,140){\makebox(0,0){$P_2$}}
\put(60,170){\makebox(0,0){$P_1$}}
\put(90,100){\makebox(0,0){$Q_{n_v}$}}
\put(90,140){\makebox(0,0){$Q_2$}}
\put(90,170){\makebox(0,0){$Q_1$}}
\put(70,50){\makebox(0,0){$G$}}
\multiput(50,90)(3,0){10}{.}
\multiput(50,130)(3,0){10}{.}
\multiput(50,160)(3,0){10}{.}

\put(40,180){\makebox(0,0){$\tilde{e}_1$}}
\put(110,170){\makebox(0,0){$\tilde{e}_2$}}

\put(240,180){\makebox(0,0){$\tilde{e}_1$}}
\put(310,170){\makebox(0,0){$\tilde{e}_2$}}

\put(310,115){\makebox(0,0){$\vdots$}}
\put(250,90){\circle*{4}}
\put(250,130){\circle*{4}}
\put(250,160){\circle*{4}}
\put(280,90){\circle*{4}}
\put(280,130){\circle*{4}}
\put(280,160){\circle*{4}}
\put(280,90){\line(1,0){40}}
\put(280,130){\line(1,0){40}}
\put(280,160){\line(1,0){40}}
\put(250,90){\line(-1,-1){25}}
\put(250,130){\line(-1,-1){25}}
\put(250,160){\line(-1,-1){25}}
\put(250,90){\line(-1,1){25}}
\put(250,130){\line(-1,1){25}}
\put(250,160){\line(-1,1){25}}
\put(260,80){\makebox(0,0){$P_{n_v}$}}
\put(260,140){\makebox(0,0){$P_2$}}
\put(260,170){\makebox(0,0){$P_1$}}
\put(290,100){\makebox(0,0){$Q_{n_v}$}}
\put(290,140){\makebox(0,0){$Q_2$}}
\put(290,170){\makebox(0,0){$Q_1$}}
\put(280,50){\makebox(0,0){$g_{e_v}*G$}}
\multiput(250,90)(2,5){15}{.}
\multiput(250,130)(3,-3){8}{.}
\multiput(250,160)(3,-3){10}{.}
\multiput(280,90)(-3,3){5}{.}

\end{picture}

Notice that there is a canonical equivariant identification of the edges
of $G$ and $g_{e_v}*G$ (which, in general, cannot be extended to a morphism of 
graphs). 
Indeed, there are natural morphisms $q_1:\tilde{G}_v\to G$ and
$q_2:\tilde{G}_v\to g_{e_v}*G$ which correspond to the gluings described above.
They induce isomorphisms at the level of edges, hence an isomorphism
of sets:
$$i_{g_{e_v}}(G):{\cale}(G)\to {\cale}(g_{e_v}*G).$$

On the above diagram of graphs, we decorated by $\tilde{e}_1$ and
$\tilde{e}_2$ two  identified pairs of edges; this partial
identification can be extended to a unique  equivariant identification
of edges.\\

It is not difficult to see, that this construction 
$(g_{e_v},G)\mapsto g_{e_v}*G$
can be extended to a group action
 of $\bgu$. Moreover, for any $g\in B(\Gamma, {\bf n})$,
 one can define inductively an identification isomorphism:
$$i_g(G): {\cale}(G)\to {\cale}(g*G).$$

\noindent 
The verification of the following easy facts are left to the reader.

\bekezdes{Lemma.}{4.8}\ {\em 
a.) If $G_1\sim G_2$ then $g*G_1\sim g*G_2$ for any $g$.
In particular, we obtain a group action
$\bgu\times\fedog\rightarrow\fedog$.

b.) The above action on $\fedog$ is transitive.}

\vspace{2mm}

\noindent
For a fixed $G$ let $Iso  (G)$ be the isotropy subgroup
$\{ g\in\bgu : g*G=G\}$ of $G$. Since $\bgu$ is abelian and the action is transitive,
all the isotropy subgroups coincide. We denote this subgroup by $Iso \subset\bgu$.
An immediate consequence of the previous lemma is that $\calg$ can be represented as a
homogeneous set $\bgu /Iso $.

In the next lemma, we will give a (possible) generator set of $Iso $.
Fix a covering $p:G\rightarrow\Gamma$. If $e=(v,w)$ is an edge of $\Gamma$, then we
can detach both endpoints of $e$ from $\Gamma$. We make the same twist above the
separated endpoints as before (but now above both endpoints).
This operation does not change the equivalence class of $G$, since we just rotate
the block $p^{-1}(e)$, which is equivariant. In other words, for any edge $e=(v,w)$,
the equivariance of block $p^{-1}(e)$ gives that $\hat{g}_e *G\sim G$, where
$\hat{g}_e = g_{e_v}\cdot g_{e_w}$.

Similarly, if $v$ is a vertex, then we take the product 
$\hat{g}_v=\prod\limits_{e_v\in{\cale}_v} g_{e_v}$ and equivariance of the fiber above
$v$ gives $\hat{g}_v * G\sim G$.

\bekezdes{Lemma.}{4.9} \ {\em 
$Iso $ is generated by the elements $\{ \hat{g}_e \}_{e\in{\cale}(\Gamma)}$ and
$\{ \hat{g}_v \}_{v\in{\calv}(\Gamma)}$.}

\vspace{2mm}

\noindent {\em Proof.} We will prove that 
$Iso(T) $ is generated by the elements $\{ \hat{g}_e \}_{e\in{\cale}(\Gamma)}$ and
$\{ \hat{g}_v \}_{v\in{\calv}(\Gamma)}$.

Fix a section $s:\G\to T$ of $T$, which determines a set of distinguished
edges $\{s(e)\}_{e\in\cale(\G)}$ of $T$ (cf. \ref{4.6}). Recall that for any 
$g\in B(\G,{\bf n})$, as we already mentioned, there is a canonical equivariant
identification  $i_g=i_g(T):\cale(T)\to \cale(g*T)$. 

Fix an element $g\in Iso(T)$, and set $s_g(e):= i_g(s(e))$. The 
correspondence $e\mapsto s_g(e)$, in general, cannot be extended to a 
section of $g*T\to \G$, (because, for a fixed $v\in\calv(\G)$, the set of 
endpoints above $v$ of 
the edges $\{s_g(e_v)\}_{e_v\in\cale_v}$ might contain more than one
vertex). 

On the other hand, since $g*T$ is trivial, there is a section $s':\G\to g*T$,
which again determines a set of distinguished edges $\{s'(e)\}_{e\in\cale(\G)}$.

For any block of $g*T$ above the segment $S_e$ of $\G$, we can find
an integer $k_e\in\bfz$ such that $e_{\cale}(k_e,s_g(e))=s'(e)$. Set 
$h=\prod_e\hat{g}_e^{-k_e}$; then the endpoints
 of the edges $\{i_{h\circ g}(s(e))\}$
above any vertex $v$ consists of exactly one vertex, i.e. they form a section
of $(h\circ g)*T$. But, the elements $g'$ of $B(\G,{\bf n})$, which
have the property that for a section $s$ of $T$, the set of edges 
$\{i_{g'}(s(e))\}_e$ can be extended to a section, have the form 
$g'=\prod_v\hat{g}_v^{k_v}$ (for some integers $k_v$). 
Hence $g=\prod \hat{g}_e^{k_e}\cdot \prod \hat{g}_v^{k_v}$.
\hfill $\diamondsuit$

\vspace{2mm}

Now, for any edge $e$ with endpoints $v_1$ and $v_2$ define  
$n_e=d_e[n_{v_1},n_{v_2}]$.
Let $\agn$ be the group 
$\prod\limits_{e\in\cale(\Gamma)}\bfz_{n_e}\times\prod
\limits_{v\in\calv(\Gamma)}{\bfz}_{n_v}$.
The generators of 
$\agn$ are $\tilde{g}_e=(0,...,\hat{1},...,0)$ (resp. $\tilde{g}_v$),
whose all entries are zero except at the place $e$ (resp. $v$), where we
put the generator $\hat{1}$  of ${\bfz}_{n_e}$ (resp. of  ${\bfz}_{n_v}$) .
Consider the map $\theta :\agn\rightarrow\bgu$ given by
$\theta(\tilde{g}_{e}) = \hat{g}_{e} = g_{e_{v_1}}\cdot g_{e_{v_2}}$, and 
$\theta(\tilde{g}_{v}) = \hat{g}_v = \prod\limits_{e_v\in\cale_v} g_{e_v}$.
Then by the above discussion $\fedog$ can be identified 
with the {\em factor group} $\coker\theta$. Indeed, 
just notice that we have 
a distinguished element of $\calg$, namely the trivial covering
$T$. The action $g\in\bgu\rightarrow g*T\in{\calg}$ identifies $\fedog$ 
with the factor  group $\coker\theta=\bgu /Iso $.

The above discussions are summarized in the following theorem:

\bekezdes{ Theorem.}{4.10}\ {\em 
a.) One has the following exact sequence of abelian groups:

$$
\agn\stackrel{\theta}{\rightarrow}\bgu\stackrel{*T}
{\longrightarrow}\fedog\rightarrow 0
$$

\noindent In particular, $\fedog$ is independent of $\bd$ 
(therefore, sometimes we will use the more simple notation 
${\calg}(\Gamma,\bn)$ for it).

b.) If $\Gamma = \coprod_{i=1}^k \Gamma _i$ has $k$ connected components
$\{\Gamma_i\} _{i=1}^k$, then

$$ \fedog = \oplus_{i=1}^k{\calg}(\Gamma_i,(\bn,\bd)).$$
}

\bekezdes{ -- ${{\calg}(\Gamma ', (\bn,\bd))}$ for a subgraph
 $\Gamma '\subset\Gamma$.}{4.13}

Fix a graph $\Gamma$ and integers $(\bn,\bd )=\{ \{ n_v\} _{v\in\calv};
\{ d_e\} _{e\in\cale}\}$ as above.
Then any subgraph $\Gamma '$ (i.e. graph  $\Gamma '$ with $\calv (\Gamma')\subset
\cvg$ and $\cale (\Gamma ')\subset\evg$ ) inherits a system of integers
$\{ \{ n_v\} _{v\in\calv (\Gamma ')};
\{ d_e\} _{e\in\cale (\Gamma ')}\}$ still denoted by $(\bn, \bd)$. The natural
projections $pr_A:\agn\rightarrow A(\Gamma ',\bn)$ and  $pr_B:\bgu\rightarrow B(\Gamma ',\bn)$ can be inserted in the following diagram:

\vspace{2cm}

\bekezdes{}{4.14}\

\

\vspace{-3.5cm}

\ 

\hspace{-0.7cm}
\begin{picture}(400,200)(30,0)
\put(130,50){\vector(1,0){40}}
\put(130,100){\vector(1,0){40}}
\put(130,150){\vector(1,0){40}}

\put(230,50){\vector(1,0){40}}
\put(230,100){\vector(1,0){40}}
\put(230,150){\vector(1,0){40}}
\put(230,50){\vector(1,0){38}}
\put(230,100){\vector(1,0){38}}
\put(230,150){\vector(1,0){38}}

\put(330,50){\vector(1,0){20}}
\put(330,100){\vector(1,0){20}}
\put(330,150){\vector(1,0){20}}

\put(100,180){\vector(0,-1){15}}
\put(100,135){\vector(0,-1){20}}
\put(100,85){\vector(0,-1){20}}
\put(100,85){\vector(0,-1){18}}
\put(100,35){\vector(0,-1){15}}

\put(200,180){\vector(0,-1){15}}
\put(200,135){\vector(0,-1){20}}
\put(200,85){\vector(0,-1){20}}
\put(200,85){\vector(0,-1){18}}
\put(200,35){\vector(0,-1){15}}

\put(300,135){\vector(0,-1){20}}
\put(300,85){\vector(0,-1){20}}
\put(300,85){\vector(0,-1){18}}
\put(300,35){\vector(0,-1){15}}

\put(100,190){\makebox(0,0){$0$}}
\put(200,190){\makebox(0,0){$0$}}
\put(100,10){\makebox(0,0){$0$}}
\put(200,10){\makebox(0,0){$0$}}
\put(300,10){\makebox(0,0){$0$}}
\put(370,150){\makebox(0,0){$0$}}
\put(370,100){\makebox(0,0){$0$}}
\put(370,50){\makebox(0,0){$0$}}

\put(100,100){\makebox(0,0){$A(\G,{\bf n})$}}
\put(100,50){\makebox(0,0){$A(\G',{\bf n})$}}
\put(200,100){\makebox(0,0){$B(\G,{\bf n})$}}
\put(200,50){\makebox(0,0){$B(\G',{\bf n})$}}
\put(300,100){\makebox(0,0){${\calg}(\G,{\bf n})$}}
\put(300,50){\makebox(0,0){${\calg}(\G',{\bf n})$}}

\put(100,150){\makebox(0,0){$\ker(pr_A)$}}
\put(200,150){\makebox(0,0){$\ker(pr_B)$}}
\put(300,150){\makebox(0,0){coker\ $\theta$}}

\put(150,160){\makebox(0,0){$\theta$}}
\put(150,110){\makebox(0,0){$\theta(\G)$}}
\put(150,60){\makebox(0,0){$\theta(\G')$}}
\put(250,60){\makebox(0,0){$*T$}}
\put(250,110){\makebox(0,0){$*T$}}

\put(115,75){\makebox(0,0){$pr_A$}}
\put(215,75){\makebox(0,0){$pr_B$}}
\put(310,75){\makebox(0,0){$pr$}}
\put(310,125){\makebox(0,0){$j$}}

\end{picture}

\noindent 
Then the horizontal lines and vertical columns are exact.
This diagram will be used several times in the next paragraphs. 
As a first immediate fact, one has:

\bekezdes{ Corollary.}{4.15}\
{\em For any subgroup $\Gamma '\subset\Gamma, \ \ 
pr:\fedog\rightarrow {\calg}(\Gamma ', (\bn,\bd))$ is onto.}

\vspace{2mm}

The construction of the coverings $G\rightarrow\Gamma$ suggests that the 
group ${\calg}(\Gamma)$ describes the possible twists above $\Gamma$.
Intuitively it is clear that a non--trivial twist
can  be realized only above cycles of $\Gamma$.
The next theorems state this fact rigorously.

\bekezdes{ Theorem.}{4.16}\ {\em 
Assume that $\Gamma$ is a tree. Then  $\fedog =0$ for any $(\bn ,\bd )$.
}

\vspace{2mm}

\noindent {\em Proof.}\ We can assume that $\G$ is connected (cf. \ref{4.10}).
We prove the vanishing of ${\calg}$
 by induction on the  number of edges of $\G$.
If $\#\cale=1$, then ${\calg}=0$ by (\ref{4.5}a) 
or by (\ref{4.10}). If $\#\cale\geq 2$,
fix a connected
subgroup $\G'\subset \G$ such that $\calv(\G')=\calv(\G)\setminus \{v\}$,
and $\cale(\G')=\cale(\G)\setminus\{(v,w)\}$
(i.e. $\Gamma'$ is obtained from $\Gamma$ by deleting a vertex
$v$ with $\#{\cale}_v=1$).  Then, in (\ref{4.14}), 
$\ker(pr_A)=
\bfz_{n_v}\times\bfz_{d_{(v,w)}[n_v,n_w]},\ \ker(pr_B)=\bfz_{n_v}\times\bfz_{n_w}$, 
and $\theta:\ker(pr_A)\to \ker(pr_B)$ is given by $(\hat{a},\hat{b})\mapsto
(\widehat{a+b},\hat{b})$, hence $\theta $ is onto. Therefore $pr:
{\calg}(\G)\to{\calg}(\G')$ is an isomorphism. Hence the theorem follows.
\hfill $\diamondsuit$


In Part II (\ref{unique}), 
 we need the following generalization of (\ref{4.16}). 

\bekezdes{Theorem.}{theo2} {\em 
Fix a graph  $\G$ and ${\bf (n,d)}$ as above.
Denote by ${\calv}^1$ the set of vertices $\{v\in {\calv}(\G):n_v=1\}$.
Let $\G_{\not=1}$ be the subgraph of $\G$ obtained from $\G$ by deleting the 
vertices ${\calv}^1$ and all those edges which have at least one endpoint in
${\calv}^1$. If  each connected component of 
$\G_{\not=1}$ is a tree, then ${\calg}(\G,{\bf (n,d)})=0$.}

\vspace{1mm}

\noindent The proof is similar to the proof of \ref{4.16} (using the 
diagram \ref{4.14}), and it is left to the reader.

\vspace{1mm}

\noindent 
Now, let $\G$ be a cyclic graph, i.e. $\calv(\G)=\{v_1,v_2,\ldots,v_k\}$,
\ $\cale(\G)=\{(v_1,v_2),(v_2,v_3),$ $\ldots,(v_k,v_1)\}$.

\bekezdes{ Theorem.}{4.17}\ {\em If \  $\G$ is a cyclic graph as above,
 and $d=g.c.d.\{n_v:
v\in\calv(\G)\}$, then ${\calg}(\G,({\bf n},{\bf d}))=\bfz_d$.}

\vspace{2mm}

\noindent 
{\em Proof.}\ Let $\G'$ be the subgraph of $\G$ with $\calv(\G')=\calv(\G)$
and $\cale(\G')=\cale(\G)\setminus\{(v_k,v_{k+1})\}$. Then in (\ref{4.14})
$\ker(pr_A)=\bfz_{d_k[n_{v_k},n_{v_{k+1}}]},\ \ker(pr_B)=\bfz_{n_{v_k}}\times
\bfz_{n_{v_{k+1}}}$, and ${\calg}(\G',{\bf n})=0$ (because $\G'$ is a tree).
Hence $\coker\theta=\bfz_{(n_{v_k},n_{v_{k+1}})}\to {\calg}(\G,{\bf n})$
is onto for any $k$. Therefore, ${\calg}(\G,{\bf n})$ is cyclic of order $d'$
with $d'|d$. 
Now, consider the map $B(\G,{\bf n})=\prod_k\bfz^2_{n_{v_k}}\to \bfz_d$
given by $((a_1,b_1),(a_2,b_2),$ $\ldots)\mapsto\sum_i a_i-\sum_ib_i$. 
This is obviously onto. The  subgroup
$Iso$ is generated by $(\ldots,(0,0),(1,1),$ $(0,0),$ $\ldots)$ and $
(\ldots,(0,1),(1,0),\ldots)$. Therefore, there is a well--defined  map 
${\calg}(\G,{\bf n})\to \bfz_d$, which is onto. Hence $d=d'$.\hfill  
$\diamondsuit$

\vspace{2mm}

\noindent The proofs of (\ref{4.16}) and (\ref{4.17}) combined give:

\bekezdes{ Corollary.}{4.18}\ {\em Let $\G$ be a graph with $c_{\G}=1$
and fix the system of integers ${\bf n}=\{n_v\}_v$. 
Let $\G'$ be the  (unique) minimal cyclic subgraph of $\G$, and set
$d:=g.c.d.\{n_v:v\in\calv(\G')\}$. Then ${\calg}(\G,{\bf n})=\bfz_d$.}

\vspace{2mm}

\noindent We end this sphere of thought with the following fact:

\bekezdes{Lemma.}{cgcg} {\em If $p:G\to \Gamma$ is a cyclic covering  then 
$c_G\geq c_{\Gamma}$.}\\
{\em Proof.} \ If $|G|$ and $|\Gamma|$ denote the topological realizations of 
the corresponding graphs, then $H_1(|G|,{\bfq})^{{\bfz}}=H_1(|\Gamma|,
{\bfq})$, hence the result follows. 
Cf. also with (\ref{7.6}), and  also with (\ref{3.4}), (\ref{9.8}) and
(\ref{9.12}). \hfill $\diamondsuit$

\bekezdes{ -- The $mod_N$--construction.}{4.19} 
(This will be crucial in Section  \ref{sec3}.)\\  
Fix an arbitrary graph $\G$ and system of integers 
${\bf n}=\{n_v\}_{v\in\calv}$ and ${\bf d}=\{d_e\}_{e\in\cale}$ as above.
For any strict positive integer $N$, we introduce  a new set of integers 
$(N,{\bf n}):=\{(N,n_v)\}_{v\in\calv}$ and ${\bf d'}$, where
$d'_e:=(d_e,N/(N,[n_{v_1},n_{v_2}]))$.
(In other words: if $n_e=d_e\cdot [n_{v_1},n_{v_2}]$, then
$(n_e,N)=d_e'\cdot [(n_{v_1},N),(n_{v_2},N)]$.) We define a natural map:
$$mod_N:{\calg}(\G,({\bf n},{\bf d}))\to {\calg}(\G,((N,{\bf n}),{\bf d'}))$$
as follows. Let $G$ be a representative of an element of 
${\calg}(\G,({\bf n},{\bf d}))$. Then 
$mod_N(G)$ is the ``orbit graph'' of the induced action of $N\bfz$.
More precisely,
we introduce the equivalence relation $\stackrel{N}{\sim}$ on $\calv(G)$, 
respectively on
$\cale(G)$: $v_1\stackrel{N}{\sim} v_2$ (resp. $e_1\stackrel{N}{\sim} e_2$) 
if there is an integer $k$
such that $a_{\calv}(kN,v_1)=v_2$ (resp. $a_{\cale}(kN,e_1)=e_2$). 
Then $mod_N(G)=G/\stackrel{N}{\sim}$.

\bekezdes{ -- Variations on the theme of coverings.}{4.20}
We can extend our set of coverings if we change the definition of the 
``standard block''. Actually, what is really important in the definition
of a block $B$ is summarized in the following two principles:

(a)\ $B$ must be equivariant (i.e. it has a ${\bfz}$--action);

(b)\ $B$ must be rigid in the following sense: if $a_{\cale}(h,e)=e$ for some 
$h\in {\bfz}$ and $e\in\cale(B)$, then all the edges and vertices of $B$ are 
left invariant by the action of $h$.

Otherwise, the block can be as complicated as we want. \\
\noindent {\bf (1)\ First variation.}\ Assume that the vertices of all our 
graphs have two types: arrowhead vertices $\cala$ and non--arrowhead 
vertices $\calw$, i.e. $\calv=\cala\coprod\calw$. Then in the definition
of the coverings $p:G\to\G$ we add the following axiom: 
\nopagebreak{$\cala(G)=p^{-1}
(\cala(\G))$.}

\noindent {\bf (2)\ Second variation.}\ Assume that our graphs have some 
decorations. Then for coverings $p:G\to\G$ we require additionally, that the
decoration of $G$ must be equivariant. (At this moment, 
we require no connection between the decorations of $G$ and $\G$.)\\
\noindent {\bf (3)\ Third variation.}\ 
The following construction is used extensively in 
Section \ref{sec3}: 
we  change (in an equivariant way)
each edge of $G$ into a string. More precisely, we start with the following
data:

(a)\ a graph $\G$ and a system ${\bf (n,d)}$ as above;

(b)\ a covering $p:G\to \G$, as an element of ${\calg}(\G,{\bf (n,d)})$;

(c)\ for each edge $e$ of $\G$, we fix a string $Str(e)$ (which may have 
some decorations):

\begin{picture}(400,50)(30,0)
\put(70,25){\makebox(0,0)[r]{$Str(e):$}}
\put(150,35){\makebox(0,0){$-k_1$}}
\put(200,35){\makebox(0,0){$-k_2$}}
\put(300,35){\makebox(0,0){$-k_s$}}
\put(150,25){\circle*{4}}
\put(200,25){\circle*{4}}
\put(300,25){\circle*{4}}
\put(225,25){\vector(-1,0){120}}
\put(275,25){\vector(1,0){70}}
\put(250,25){\makebox(0,0){$\cdots$}}
\end{picture}

Denote the collection $\{Str(e)\}_{e\in\cale(\G)}$ of strings  
by ${\bf Str}$.

\noindent
Then the new graph $G({\bf Str})$ is constructed as follows: 
we replace each 
edge $\tilde{e}\in p_{\cale}^{-1}(e)$
(with endpoints $\tilde{v}_1$ and $\tilde{v}_2$) of $G$  with the 
decorated string
$\overline{Str}(e)$, as it is shown below:

\hspace{-0.5cm}
\begin{picture}(400,50)(0,0)
\put(0,25){\circle*{3}}
\put(30,25){\circle*{3}}
\put(0,15){\makebox(0,0){$\tilde{v}_1$}}
\put(30,15){\makebox(0,0){$\tilde{v}_2$}}
\put(50,35){\makebox(0,0)[l]{\mbox{replaced}}}
\put(60,20){\makebox(0,0)[l]{\mbox{by}}}
\put(150,35){\makebox(0,0){$-k_1$}}
\put(200,35){\makebox(0,0){$-k_2$}}
\put(300,35){\makebox(0,0){$-k_s$}}
\put(100,15){\makebox(0,0){$\tilde{v}_1$}}
\put(350,15){\makebox(0,0){$\tilde{v}_2$}}

\put(100,25){\circle*{4}}
\put(150,25){\circle*{4}}
\put(350,25){\circle*{4}}
\put(200,25){\circle*{4}}
\put(300,25){\circle*{4}}
\put(100,25){\line(1,0){120}}
\put(350,25){\line(-1,0){70}}
\put(0,25){\line(1,0){30}}

\put(250,25){\makebox(0,0){$\cdots$}}
\end{picture}

\section{Graphs associated with analytic  germs}\label{sec2}

Let $(X,x)$ be a normal surface singularity and fix the germ
$f:(X,x)\to (\bfc,0)$ of an analytic function. In the first subsection we 
review the definition of the embedded resolution graph $\Gamma(X,f)$ of $f$,
and we introduce  our basic notations. 
For more details the reader is invited to consult the 
books of H. Laufer \cite{La} and Eisenbud--Neumann \cite{EN}, as well as
the survey article of J. Lipman \cite{Lip}.
In the next subsection, we
discuss Jung's method of desingularization of normal surface singularities.
Here we also recall the basic arithmetical properties of the resolution graph
of Hirzebruch--Jung singularities. The third subsection
deals with the topology of the 
link $L_X$ of $(X,x)$ and of the pair $(X,f^{-1}(x))$. We discuss in details 
the representation $arg_*(f)$ provided by the Milnor fibration associated with
$f$. This representation will be crucial in the next subsection, and in
Section  \ref{sec3}. 
In the last subsection we introduce the ``universal cyclic 
$\bfz$--covering'' $G(X,f)$ of the embedded resolution $\Gamma(X,f)$. 
This covering graph coordinates all the resolution graphs of the cyclic
coverings of $(X,x)$ branched along $f^{-1}(0)$ (see the next section).

In the body of this section we present many examples in order to clarify the 
relationship between the objects $\Gamma(X,f)$, $arg_*(f)$, $L_X$ and $G(X,f)$.
The verification of these examples sometimes is not absolutely trivial,
but basically all of them are based on the main result (\ref{7.2}).
Even if that theorem is presented only in the next section, we 
prefer to give all these examples   already here. 
The reader can skip the verification of them  at the first reading.

\subsection*{The embedded resolution graph $\G(X,f)$.}

\vspace{2mm}

\bekezdes{-- The embedded resolution.}{1.1} 
Let $(X,x)$ be a normal surface singularity and 
let $f:(X,x)\to (\bfc,0)$ be the germ of an analytic function.

An embedded resolution  $\phi:({\caly},D)\to (U,f^{-1}(0))$ of
$(f^{-1}(0),x)\subset (X,x)$  is characterized
by the following properties.  There is a sufficiently small neighborhood 
$U$ of $x$ in $X$, smooth analytic manifold ${\caly}$, and an analytic proper
map $\phi:{\caly}\to U$ such that:

1)\ if $E=\phi^{-1}(x)$, then the restriction $\phi|_{{\caly}\setminus E}:
{\caly}\setminus E\to U\setminus\{x\}$ is biholomorphic, and ${\caly}\setminus E$ is dense in ${\caly}$;

2)\ $D=\phi^{-1}(f^{-1}(0))$ 
is a divisor with only normal crossing singularities, i.e. 
at any point $P$ of $E$, there are local coordinates $(u,v)$ in some small
neighbourhood of $P$, such that in these coordinates
$f\circ\phi=u^av^b$ for some non--negative integers $a$ and $b$.

If such an embedded resolution $\phi$ is fixed, then 
$E=\phi^{-1}(x)$ is called the  exceptional divisor associated  with $\phi$.
Let $E=\cup_{w\in \calw}E_w$ be its decomposition  in irreducible divisors.
The closure $S$ of $\phi^{-1}(f^{-1}(0)\setminus\{0\})$ is called the 
strict transform of $f^{-1}(0)$.
Let $\cup_{a\in \cala}S_a$ be its irreducible decomposition.
Obviously,  $D=E\cup S$. 

In this section, for simplicity,  we will assume that $\calw\not=\emptyset$, 
any two irreducible components of $E$ have at most one intersection point, 
and  no irreducible exceptional divisor  has a self--intersection.
This can always be realized by some additional blow-ups.

\bekezdes{-- The embedded resolution graph.}{1.2}
We construct the embedded resolution graph  $\G(X,f)$ of the pair $(X,f)$,
associated with a fixed  resolution $\phi$,  as follows. Its vertices
$\calv=\calw\coprod\cala$
 consist of the nonarrowhead vertices $\calw$  corresponding to the
irreducible exceptional divisors, and arrowhead vertices $\cala$
correponding to the irreducible divisors of the strict transform $S$.
If two irreducible divisors corresponding to $v_1,v_2\in \calv$
have an intersection point then $(v_1,v_2)\ (=(v_2,v_1))$ is an edge of
$\G(X,f)$. The set of edges is denoted by $\cale$. 

For any $w\in \calw$, we denote by $\calv_w$ the set of vertices
$v\in \calv$ adjacent to $w$. Its cardinality $\#\calv_w$ is called the 
degree $\delta_w$ of $w$. In Section \ref{sec1}  
 we introduced for any graph
the set ${\cale}_w$ of 
edges adjacent to $w$. Since, by our assumption,
 any two vertices are connected by at most one edge, one has the 
identification ${\cale}_w={\calv}_w$. 

The graph $\G(X,f)$ is decorated as follows. Any  vertex $w\in \calw$ is 
decorated with three  numbers. The first is the 
self--intersection   $e_w:=E_w\cdot E_w$. 
Equivalently, $e_w$ is the Euler--number of the normal bundle
of $E_w$ in ${\caly}$. The second is
the genus $g_w$ of $E_w$.
The third decoration is given by the multiplicity of $f$. More precisely, 
for any $v\in\calv$,
 let $m_v$ be the vanishing order of $f\circ \phi$ along the
irreducible divisor corresponding to $v$. For example, if $f$ defines
an isolated singularity, then  for any
$a\in\cala$ one has $m_a=1$.

In all our graphs, we put the multiplicities in parentheses
(e.g.: (3)) and the genera in brackets (e.g.: [3]).
In order to simplify the graphs, if $g_w=0$ for some $w$, then we omit $[0]$
from the graph.

\bekezdes{-- First properties of the graph.}{1.3}\ 

(1)\ Notice that the combinatorics of the 
graph and the self-intersection numbers
codify completely the intersection matrix $(E_w\cdot E_v)_{(w,v)\in\calw\times
\calw}$ of the irreducible components of $E$. Moreover, this matrix is 
negative definite, see e.g. 
\cite{Mumford} page 230; \cite{La} page 49, or 
\cite{GRa}.

(2)\ The following compatibility  property holds between self--intersections 
and multiplicities. For any $w\in\calw$ one has:
$$e_wm_w+\sum_{v\in\calv_w}m_v=0.$$

Obviously, $m_v>0$ for any $v\in{\calv}$, 
hence the set of multiplicities determine the 
self--intersection numbers completely. 
The advantage of this fact is the following: a multiplicity can always be determined by a local computation, on the other hand the Euler--number $e_w$ is a
global characteristic class.

Similarly, since the intersection 
matrix $(E_w\cdot E_v)_{(w,v)\in\calw\times\calw}$
is negative definite,  these  relations  determine the 
multiplicities $\{m_w\}_{w\in\calw}$ in terms
of the self--intersection numbers  $\{e_w\}_w$
and the multiplicities $\{m_a\}_{a\in\cala}$.

\bekezdes{-- The resolution of $(X,x)$.}{1.4}

We say that $\phi:{\caly}\to U$ is a resolution of $(X,x)$ if 
${\caly}$ is a smooth analytic manifold, $U$ a neighbourhood of $x$ in $X$, 
$\phi$ is a proper analytic map, such that 
${\caly}\setminus E$ (where $E=\phi^{-1}(x)$) is dense in ${\caly}$ and
the restriction $\phi|_{{\caly}\setminus E}:{\caly}\setminus E\to U\setminus
\{x\}$ is biholomorphic.

The topology  of the resolution and the combinatorics of the  irreducible
exceptional divisors $\cup_wE_w$ 
can by codified in the graph $\G(X)$, which is called the dual 
resolution graph of $(X,x)$ associated with $\phi$. 
If the divisor $E$ is not a normal crossing divisor, then this codification 
can be slightly complicated, so in the sequel we will assume that the 
irreducible components of $E$ are smooth and intersect each other 
transversally, the irreducible components have no self--intersections,
and there is no intersection point which is contained in 
more than two components.
In this case, similarly as in the situation of the embedded resolution,
the vertices of the dual graph
correspond to the irreducible components of $E$, the edges 
to the intersections of these components, and each vertex $w$ carry two
decorations: the genus $g_w$ of $E_w$, and the self--intersection
$E_w^2$. 

Actually, one can obtain 
a possible graph $\G(X)$ from any $\G(X,f)$ by deleting 
all the arrows and multiplicities of the graph $\G(X,f)$.

The graphs $\G(X,f)$ and $\G(X)$ are connected (see  \cite{La}, 
or Zariski's Main Theorem, e.g. in \cite{Hartshorne}).

\bekezdes{Example. -- Plane curve singularities.}{1.6}
(see e.g. \cite{BrKn}) \
Assume that $(X,x)$ is smooth. Then the singular germ
$(f^{-1}(x),x)\subset (X,x)$ can be resolved only by quadratic modifications. 
In this case, the graph $\G(X,f)$ is a tree, and $g_w=0$ for
any $w\in\calw$. 

\bekezdes{ Example. -- The normalization of an arbitrary surface singularity.}{1.7}
If $(X,x)$ is a surface singularity, but it is not normal, then there is a
canonical way to construct its normalization $\tilde{X}$ (see e.g. \cite{La}).
If $(X,x)$ has $k$ local irreducible components, then its 
normalization will split into $k$ normal singular space--germs. 
The corresponding 
resolution graph of $(\tilde{X},\{x_1,\ldots,x_k\})$ 
will have $k$ connected components, each 
corresponding exactly to a resolution graph of the irreducible components
$(\tilde{X},x_i)$. 

\bekezdes{ Example. -- Cyclic coverings.}{1.10}
Start with a normal surface singularity $(X,x)$ and a germ $f:(X,x)\to 
({\bfc},0)$. Consider the covering $b:({\bfc},0)\to ({\bfc},0)$ given by
$z\mapsto z^N$, and construct the fiber product:
$$(X,x)\prod_{f,b}({\bfc},0)=\{(x',z)\in(X\times{\bfc},x\times 0)\ :\ f(x')=
z^N\}.$$
By definition, $X_{f,N}$ is the normalization of 
$(X,x)\prod_{f,b}({\bfc},0)$.
There is a natural  ramified covering $X_{f,N}\to X$ branched along
$f^{-1}(0)$.

\bekezdes{ Example. -- Hirzebruch--Jung singularities.}{1.9}
(See \cite{Hir1,La,BPV}). 
For a normal surface singularity, the following conditions are equivalent.
If $(X,x)$ satisfies (one of) them, then it is called Hirzebruch--Jung 
singularity.

a)\ The resolution graph $\G(X)$ is a string, and  $g_w=0$ for any
$w\in\calw$. (If the graph is minimal then 
 $e_w\leq -2$ for any $w$.)

b)\ There is a finite proper map
$\pi:(X,x)\to ({\bfc}^2,0)$ such that the reduced discriminant locus 
of $\pi$, in some local coordinates $(u,v)$ of $({\bfc}^2,0)$, is $\{uv=0\}$.

c)\ $(X,x)$ is isomorphic with exactly one of the ``model spaces'' 
$\{A_{n,q}\}_{n,q}$, where $A_{n,q}$ is the normalization of 
 $(\{xy^{n-q}+z^n=0\},0)$, with $0<q<n,\ (n,q)=1$.\\

If there is a map $\pi$ as in (b) with smooth reduced discriminant locus, then 
$(X,x)$ is automatically smooth.

\subsection*{Jung's method. Hirzebruch--Jung singularities.}

\vspace{2mm}

\bekezdes{}{2.1}
Jung's method \cite{Jung} gives not only a qualitative
proof of the existence of the resolution $\phi:({\caly},E)\to (X,x)$ of a 
normal surface singularity $(X,x)$ (see, e.g. \cite{La} Theorem 2.1;  or
\cite{BPV} Theorem 6.1), but also  a rather clear recipe how this resolution
 can be  constructed in concrete cases (see, e.g. \cite{La} chapter III,
 or \cite{Lip}). 
The Jungian strategy (presented for the case of a normal surface 
singularity $(X,x)$) can be summarized in the following diagram:

\begin{picture}(400,100)(0,0)
\put(50,75){\makebox(0,0){$X^{res}$}}
\put(100,75){\makebox(0,0){$\tilde{X}'$}}
\put(160,75){\makebox(0,0){$X'$}}
\put(230,75){\makebox(0,0)[l]{$(X,x)$}}
\put(160,25){\makebox(0,0){$({\caly},D)$}}
\put(230,25){\makebox(0,0)[l]{$(\bfc^2,0)\supset (\Delta,0)$}}
\put(65,75){\vector(1,0){20}}
\put(115,75){\vector(1,0){20}}
\put(190,75){\vector(1,0){30}}
\put(190,25){\vector(1,0){30}}
\put(115,60){\vector(1,-1){20}}
\put(160,60){\vector(0,-1){20}}
\put(250,60){\vector(0,-1){20}}
\put(60,60){\vector(2,-1){50}}
\put(60,50){\makebox(0,0){$\rho$}}
\put(170,50){\makebox(0,0){$\pi'$}}
\put(115,50){\makebox(0,0){$\tilde{\pi}'$}}
\put(260,50){\makebox(0,0){$\pi$}}
\put(200,85){\makebox(0,0){$\phi_{\Delta}'$}}
\put(200,35){\makebox(0,0){$\phi_{\Delta}$}}

\end{picture}

\noindent where:

a)\ $\pi:(X,x)\to (\bfc^2,0)$ is a proper finite map with (reduced)
discriminant locus $(\Delta,0)\subset (\bfc^2,0)$. 

b)\ $\phi_{\Delta}$ is an embedded resolution of $(\Delta,0)\subset 
(\bfc^2,0)$. In particular, $D=\phi_{\Delta}^{-1}(\Delta) $ is a divisor 
with normal  crossing singularities.

c)\ $\pi':X'\to {\caly}$ is the pullback of $\pi$ via $\phi_{\Delta}$, and
$\tilde{X}'$ is the normalization of $X'$. Then $\tilde{X}'$ has only normal 
singularities, and the discriminant of the projection $\tilde{\pi}'$ has 
normal crossing singularities only. This property characterizes exactly the 
Hirzebruch--Jung singularities (cf. \ref{1.9}). 

d)\ $X^{res}\to \tilde{X}'$ is the resolution of the (Hirzebruch--Jung)
singularities of $\tilde{X}'$.

\bekezdes{Example.}{2.2} In general, the computation of the 
discriminant locus $\Delta$ of some projection $\pi$, and the whole process, 
can be rather complicated. But in some cases, the above 
recipe is really nice. For example, if $(X,x)=(\{f(x,y)+z^N=0\},0)\subset
(\bfc^3,0)$, and $\pi$ is induced by $(x,y,z)\mapsto (x,y)$, then 
$(\Delta,0)=(\{f=0\},0)\subset (\bfc^2,0)$. It turns out,
that  the dual resolution graph of $(X,x)$ can be
recovered from the embedded resolution graph of $f$ and the integer $N$
(see  \ref{8.1}).

\vspace{2mm}

The above strategy shows that, in order to resolve $(X,x)$, we have to know
two things: the embedded resolution of plane curve singularities,
and the resolution of Hirzebruch--Jung singularities.

\bekezdes{-- Hirzebruch--Jung singularities.}{2.3} \cite{Hir1,La,BPV}
(cf. also \ref{1.9} ).

From our point of view, (and also from the point of view of the strategy
\ref{2.1}), it is more convenient to consider a biger class of ``models'' instead
of $\{A_{n,q}\}_{n,q}$. 

For any three strictly positive integers $a,b$ and $N$, with $g.c.d.(a,b,N)=1$,
we define $(X,x)=
(X(a,b,N),x)$ as the unique singularity lying over the origin in 
the normalization of $(\{\alpha^a\beta^b+\gamma^N=0\},0)$. 
Let the germ $\gamma:(X(a,b,N),x)\to (\bfc,0)$ be 
induced by $(\alpha,\beta,\gamma)\mapsto \gamma$. 
In the sequel, we  give the embedded resolution
graph $\G(X,\gamma)$ of the germ $\gamma$. 
Obviously, if we delete the arrows and 
multiplicities of $\G(X,\gamma)$ we obtain the resolution graph $\G(X)$ of 
$(X(a,,b,N),x)$.

First,  consider  the unique $0\leq \lambda <N/(a,N)$  and 
$m_1\in{\bf N}$ with:
$$b+\lambda \cdot \frac{a}{(a,N)}=m_1\cdot \frac{N}{(a,N)}.$$
If $\lambda\not=0$, then consider the continuous fraction:
$$\frac{N/(a,N)}{\lambda}=k_1-{1\over\displaystyle
k_2-{\strut 1\over\displaystyle\ddots
-{\strut 1\over k_s}}}, \ \ k_1,\ldots, k_s\geq 2.$$
Then the following string, denoted by $Str(a,b;N)$,
 is the embedded resolution graph $\G(X,\gamma)$ of $\gamma$:

\hspace{-2.8cm}
\begin{picture}(400,50)(-20,0)
\put(95,25){\makebox(0,0)[r]{$(\frac{a}{(a,N)})$}}
\put(355,25){\makebox(0,0)[l]{$(\frac{b}{(b,N)})$}}
\put(150,35){\makebox(0,0){$-k_1$}}
\put(200,35){\makebox(0,0){$-k_2$}}
\put(300,35){\makebox(0,0){$-k_s$}}
\put(150,15){\makebox(0,0){$(m_1)$}}
\put(200,15){\makebox(0,0){$(m_2)$}}
\put(300,15){\makebox(0,0){$(m_s)$}}
\put(150,25){\circle*{4}}
\put(200,25){\circle*{4}}
\put(300,25){\circle*{4}}
\put(225,25){\vector(-1,0){120}}
\put(275,25){\vector(1,0){70}}
\put(250,25){\makebox(0,0){$\cdots$}}
\end{picture}

\noindent 
The arrow on the left (resp. right) hand side codifies the strict
transform of $\{\alpha=0\}$ (resp. of $\{\beta=0\}$). 
All vertices have genus $g_w=0$, i.e. they represent rational
irreducible exceptional divisors. 
The first vertex has
multiplicity $m_1$ given by the above congruence. 
Hence $m_2,\ldots, m_s$
can  easily be computed using (\ref{1.3}), namely:
$$-k_1m_1+\frac{a}{(a,N)}+m_2=0,\ \mbox{and}\
-k_im_i+m_{i-1}+m_{i+1}=0\ \mbox{for $i\geq 2$}.$$

\noi 
This resolution resolves also the germ $\alpha$ (induced by  the projection
$(\alpha,\beta,\gamma)$ $\mapsto \alpha$).
The multiplicities of $\alpha$ along the (same!) divisors
(exceptional divisors
and strict transforms of $\{\alpha=0\}$ and $\{\beta=0\}$)
are given in the next graph. (This is, in fact, the
graph $\G(X,\alpha)$, if we delete the arrow with zero multiplicity):

\hspace{-0.7cm}
\begin{picture}(400,50)(30,0)
\put(355,25){\makebox(0,0)[l]{$(0)$}}
\put(150,35){\makebox(0,0){$-k_1$}}
\put(200,35){\makebox(0,0){$-k_2$}}
\put(300,35){\makebox(0,0){$-k_s$}}
\put(150,15){\makebox(0,0){$(\lambda)$}}
\put(95,25){\makebox(0,0)[r]{$(\frac{N}{(a,N)})$}}
\put(150,25){\circle*{4}}
\put(200,25){\circle*{4}}
\put(300,25){\circle*{4}}
\put(300,15){\makebox(0,0){$((b,N))$}}
\put(275,25){\vector(1,0){70}}
\put(225,25){\vector(-1,0){120}}
\put(250,25){\makebox(0,0){$\cdots$}}
\end{picture}

\noindent
The other multiplicities can  again be computed by (\ref{1.3}).
We emphasize again: the arrows codify
the {\em same strict transforms} as the arrows of $\G(X,\gamma)$.

The multiplicity of $\beta$ along the corresponding irreducible divisors can
be determined symmetrically:

\hspace{-1cm}
\begin{picture}(400,50)(30,0)
\put(95,25){\makebox(0,0)[r]{$(0)$}}
\put(355,25){\makebox(0,0)[l]{$(\frac{N}{(b,N)})$}}
\put(150,35){\makebox(0,0){$-k_1$}}
\put(200,35){\makebox(0,0){$-k_2$}}
\put(300,35){\makebox(0,0){$-k_s$}}
\put(150,15){\makebox(0,0){$((a,N))$}}
\put(150,25){\circle*{4}}
\put(200,25){\circle*{4}}
\put(300,25){\circle*{4}}
\put(300,15){\makebox(0,0){$(\tilde{\lambda})$}}
\put(225,25){\vector(-1,0){120}}
\put(275,25){\vector(1,0){70}}
\put(250,25){\makebox(0,0){$\cdots$}}
\end{picture}

\noi where $0\leq \tilde{\lambda}<N/(N,b)$ and
$$a+\tilde{\lambda}\cdot \frac{b}{(b,N)}=m_s\cdot\frac{N}{(b,N)}.$$
Obviously, the embedded resolution graph $\G(X,\alpha^i\beta^j\gamma^k)$
of the function $\alpha^i\beta^j\gamma^k$ defined on $X$ can be deduced easily
from the above resolution graphs. It has the same shape, the same 
self--intersections and genera, and the multiplicity $m_v$ (for any vertex
$v$) satisfies:
$$m_v(\alpha^i\beta^j\gamma^k)=i\cdot m_v(\alpha)+j\cdot m_v(\beta)+
k\cdot m_v(\gamma).$$

\noi {\em Notation.} The decorated string 
$\G(X,\alpha^i\beta^j\gamma^k)$  will sometimes be denoted by:
$$Str(a,b;N\,|\,i,j;k).$$

Form the point of view of the classification theorem (\ref{1.9}c),
$X(a,b,N)$ is an $A_{n,q}$--singularity,
where $n=N/(a,N)(b,N)$ and $q=\lambda/(b,N)$ 
(cf. e.g. \cite{BPV}, page 83-84).\\

If $\lambda=0$, then the string has no vertices, in particular
$(X(a,b,N),x)$ is smooth. Moreover, in this case, the zero set of $\gamma$ 
(on $X$) has only 
a normal crossing singularity: in some local coordinates $(u,v)$ of $(X,x)$,
it can be represented as $\gamma=u^{a/(a,N)}v^{b/(b,N)}$. (In this case the 
above string becomes a double arrow, without any non--arrowhead vertices.)

\subsection*{The topology of the link of $f$.}

\vspace{2mm}

\bekezdes{-- The link of $(X,x)$.}{3.1}
Let $(X,x)$ be a normal surface singularity, and fix an embedding
$(X,x)\subset ({\bfc}^N,0)$ for some $N$. Then, for sufficiently small
$\epsilon_0>0$, all the spheres $S_{\epsilon}=\{z\in\bfc^N\ :\ ||z||=
\epsilon\}$ ($0<\epsilon\leq\epsilon_0$) intersect $(X,x)$ transversally
(see, e.g. \cite{Milnorbook,Lo}), and the differentiable manifold $S_{\epsilon}
\cap X$ does not depend on the choice  of $\epsilon$ and of the embedding 
$(X,x)\subset\bfc^N,0)$. 
It  inherits  a natural orientation and it is always connected.
It is called {\em the link of  $(X,x)$}, and denoted by $L_X$. 

From a topological point of view, $L_X$ characterizes  $(X,x)$ completely.
If $B_{\epsilon}$ denotes the ball 
$\{z\in\bfc^N\ :\ ||z||\leq \epsilon\}$, then for $\epsilon$
sufficiently small $(B_{\epsilon}\cap X,x)$ is homeomorphic to 
$(\mbox{Cone}(L_X),\mbox{vertex of the cone})$. 
We will write $U=B_{\epsilon}\cap X$ for a small  $\epsilon$.
For such a  $U$, 
consider an embedded resolution $\phi:{\caly}\to U$. Then the inclusion
$\phi^{-1}(0)=E\hookrightarrow {\caly}$ admits a strong deformation retract
$r:{\caly}\to E$, and ${\caly}$ is a manifold with smooth boundary.
Moreover, the restriction of $\phi$ to $\partial {\caly}$ identifies 
$\partial {\caly}$ with $L_X=\partial U$. This shows that 
$L_X$ is  the {\em plumbed manifold} $M(L_X)$  associated with the graph
$\G(X)$ (for details, see  \cite{NePl}), i.e. 
$\G(X)$ determines completely the 3--manifold $L_X$.
The converse is also true: W. Neumann in \cite{NePl} proved that the topology
of the (minimal)  resolution of the   singularity $(X,x)$ 
is determined  by the oriented homeomorphism type
of the link $L_X$. 

\bekezdes{-- The homology of $L_X$.}{3.3}
Consider the intersection 
matrix $A$ given by $A_{v,w}=E_v\cdot E_w$. 
This defines a bilinear form $(\bfz^{\#\calw})^{\otimes 2}\to\bfz$, or
equivalently a $\bfz$--linear map $A:\bfz^{\#\calw}\to (\bfz^{\#\calw})^*$
(where for a $\bfz$--module $M$, $M^*$ denotes its dual $Hom_{\bfz}(M,\bfz)$).
Since $A$ is non--degenerate, $\coker A$ is a torsion group with 
$|\coker A|=|\det A|$. Then the following holds:

\bekezdes{Proposition.}{3.4} \cite{Sch,HNK,Mumford} \ \ 
$H_1(L_X,\bfz)\approx \coker A\oplus \bfz^{2g+c_{\Gamma}}.$

\vspace{2mm}

In particular, $L_X$ is an integer (resp. rational) homology sphere 
if and only if  $g=c_{\Gamma}=0$ and $\det A=\pm 1$ (resp. $g=c_{\Gamma}=0$;
i.e. $\Gamma$ is a tree with $g_w=0$ for all $w$).

\bekezdes{-- The links of germs $f:(X,x)\to (\bfc,0)$.}{3.6}
By similar arguments (and notations) as above, for sufficiently small
$\epsilon$ the intersection $S_{\epsilon}\cap f^{-1}(0)\subset S_{\epsilon}
\cap X=L_X$ defines a 1--dimensional compact (in general non--connected)
orientable submanifold $L_f$ of $L_X$, called the link of $f$.
The link--components (=connected components of $L_f$) are indexed by  $\cala$
in $\G(X,f)$. Similarly as above, one can recover the link 
$L_f\subset L_X$ from the graph $\G(X,f)$
by a plumbing construction.

\bekezdes{ -- The homology group $H_1(L_X\setminus L_f,\bfz)$.}{3.7}
Let $\bfz^{\calv}$ be the free abelian group generated by $\{[v]\}_{v\in\calv}$
(recall $\calv=\calw\coprod\cala$). 
Define the group $H_{\G}$ as the quotient of
$\bfz^{\calv}$ factorized by the subgroup generated by:
$$e_w[w]+\sum_{v\in\calv_w}[v]\ \ \ (\mbox{for all $w\in\calw$}).$$
Let $i$ be the composed map $\bfz^{\cala}\hookrightarrow \bfz^{\calv}\to 
H_{\G}$. Then one has the following exact sequence:

\bekezdes{}{3.8}
\hspace{2cm}$0\to\bfz^{\cala}\stackrel{i}{\to}H_{\G}\to \coker A\to 0.$

\bekezdes{Proposition.}{3.9}\ (\cite{Sch,HNK,NePl}) 
{\em There is a natural exact sequence:}
$$0\to H_{\G}\stackrel{j}{\to} H_1(L_X\setminus L_f,\bfz)\stackrel{q}{\to}
H_1(E,\bfz)\to 0.$$

Geometrically, the map $j$ can be described as follows. Identify $L_X$ with 
$\partial {\caly}$ and $L_f$ with the  intersection 
$\cup_aS_a\cap L_X$ of $L_X$ with
the strict transform $S$. For each $v\in\calv$, define 
$M_v\subset L_X\setminus L_f$ as a naturally oriented circle in the transversal
slice of the corresponding irreducible divisor of $D=
(\cup_wE_w)\cup(\cup_aS_a)$. Equivalently, for 
$w\in\calw$ we can take $M_w$ as 
a generic fiber  of $T_w\to B_w$, and for $a\in\cala$, $M_a$ is a topological 
standard meridian (see e.g. \cite{EN}) 
of $S_a\cap L_X \subset L_X$. Then $j(\widehat{[v]})
=[M_v]$ for any $v\in\calv$.

If $S$ is the strict transform, then the strong deformation retract
$r:{\caly}\to E$ can be chosen
such that it preserves $S$, hence induces a map
$r':{\caly}\setminus D\to E\setminus S$. 
Then $q$ is the composed map:

$H_1(L_X\setminus L_f,\bfz)\stackrel{\phi_*^{-1}}{\approx}
H_1({\caly}\setminus D,\bfz)\stackrel{r'_*}{\longrightarrow}
H_1(E\setminus S,\bfz)\to H_1(E,\bfz).$

\bekezdes{-- The Milnor fibration.}{3.10}\ Fix a map $f:(X,x)\to (\bfc,0)$.
Then $arg=f/|f|:L_X\setminus L_f\to S^1$ is a 
${\calc}^{\infty}$ fibration (\cite{Milnorbook,Lo}).
This induces $\pi_1(arg):\pi_1(L_X\setminus L_f)\to \bfz$ at the fundamental 
group level, and $arg_*:H_1(L_X\setminus L_f,\bfz)\to \bfz$ at homology level.
Obviously, if $ab:\pi_1\to \pi_1/[\pi_1,\pi_1]=H_1$ is the abelianization map,
 then  $\pi_1(arg)=ab\circ arg_*$. 

The connection with the exact sequence (\ref{3.9}) is the following. 
Let ${\bf m}:
\bfz^{\calv}\to \bfz$ be defined by ${\bf m}([v])=m_v$ (where $m_v$ is the 
multiplicity of $f$ along $E_w$ as in \ref{1.2}). Then by (\ref{1.3}), 
this induces a 
well--defined map ${\bf m}:H_{\G}\to\bfz$. Then $arg_*\circ j={\bf m}$
(cf. \cite{EN}):

\begin{picture}(300,75)(0,-10)

\put(150,50){\makebox(0,0){$
0\to H_{\G}\stackrel{j}{\to} H_1(L_X\setminus L_f,\bfz)\stackrel{q}{\to}H_1(E,\bfz)\to 0.$}}

\put(130,35){\vector(0,-1){20}}
\put(150,25){\makebox(0,0){$arg_*$}}
\put(85,35){\vector(2,-1){40}}
\put(80,25){\makebox(0,0){${\bf m}$}}

\put(130,0){\makebox(0,0){$\bfz$}}
\end{picture}

In particular, $arg_*$ contains all the information 
about the multiplicities, but as we will see later, from the multiplicities
we cannot recover the representation $arg_*$ (cf. Examples 
\ref{3.14}--\ref{3.18}). 
On the other hand, in the study of cyclic coverings of $(X,x)$
branched along $f^{-1}(0)$, the representation $arg_*$ plays a 
crucial role. In \ref{5.1},  
the additional information about
 $arg_*$ which is not contained in $\Gamma(X,f)$ will be codified in 
a $\bfz$--covering graph of $\Gamma(X,f)$. 
Moreover, one can verify the following:

\bekezdes{Proposition.}{3.12} {\em The fibration $arg:L_X\setminus L_f\to
S^1$ is completely  determined (up to isotopy) 
by  the induced representation $arg_*:H_1(L_X\setminus L_f,\bfz)\to \bfz$.
Moreover, if  $\bfz_d:=
\coker(p_*:H_1(L_X\setminus L_f,\bfz)\to \bfz)$, then the (Milnor)
fiber of $arg$ has $d$ 
connected components which are cyclically permuted by the monodromy.

In general, one has the following divisibility conditions:
$d\,|\,g.c.d.\{m_v:v\in\calv\}\,|\,g.c.d.\{m_a:a\in\cala\}$, but it is possible
that $d\not= g.c.d.\{m_v:v\in\calv\}$. This means that $d$ cannot be determined from ${\bf m}$.}\\

Since ${\bf m}$, in general,  does not determine the representation
$arg_*$ (cf. \ref{3.10}), the graph $\G(X,f)$ 
alone does not determine the Milnor fibration associated with $f$.

\bekezdes{-- The particular case when  
$L_X$ is a rational homology sphere.}{3.13}\
Set $f:(X,x)\to (\bfc,0)$ as above, and assume that $L_X$ is 
a rational homology sphere. This is equivalent to the vanishing of
$H_1(E,\bfz)$. Then by (\ref{3.9}\ --\ 
\ref{3.10}), $\pi_1(arg):\pi_1(L_X\setminus L_f)\to\bfz$
is completely determined by ${\bf m}:H_{\G}\to\bfz$, and  ${\bf m}$ is 
determined by $\{m_a\}_{a\in\cala}$ via 
(\ref{1.3}). Hence, the set $\{m_a\}_{a\in\cala}$
determines completely the Milnor fibration up to an isotopy. Moreover,
if $d=g.c.d.\{m_v:v\in\calv\}$, then the fiber has $d$ connected
components, and they are cyclically permuted by the monodromy.\\

The next examples show that these properties are not true if $L_X$ is not
a rational homology sphere, i.e. even with the same multilink,
different representations $H_1(L_X\setminus L_f,\bfz)\to \bfz$ do  occur.

\bekezdes{ Example.}{3.14}\ 
Set $(X,x)=(\{x^2+y^7-z^{14}=0\},0)\subset
(\bfc^3,0)$ and take $f_1(x,y,z)=z^2$ and $f_2(x,y,z)=z^2-y$. Then 
$\G(X,f_1)=\G(X,f_2)$ is :

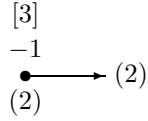
\begin{picture}(400,60)(0,0)
\put(110,20){\circle*{4}}
\put(110,30){\makebox(0,0){$-1$}}
\put(110,43){\makebox(0,0){$[3]$}}
\put(110,20){\vector(1,0){30}}
\put(110,10){\makebox(0,0){$(2)$}}
\put(150,20){\makebox(0,0){$(2)$}}
\end{picture}

Since $\coker({\bf m})=\bfz_2$, in both cases $\coker(arg_*)$ is a factor 
group of $\bfz_2$ (cf. \ref{3.10}). We will show that
in the first case $\coker(arg_*)=\bfz_2$ and in the second case $arg_*$ is 
onto. Indeed, the Milnor fibration of $z^2$ is the pullback by $z\mapsto z^2$
of the Milnor fibration of $z$, hence $\coker(arg_*)=\bfz_2$.
In order to prove the second statement, it is enough to verify that 
the double covering
$\{x^2+y^7-z^{14}=w^2+y-z^2=0\}\subset \bfc^4$ is irreducible
(notice that our equations are quasi--homogeneous, so we can replace a small
ball centered at the origin with the whole affine space). But this is 
true if its intersection with $y=1$, i.e. $C:=\{x^2=z^{14}-1;w^2=z^2-1\}
\subset \bfc^3$,  is irreducible. The covering $C\to \bfc$ ($(x,w,z)\mapsto z$)
is a $\bfz_2\times\bfz_2$ covering. The monodromy around $\pm 1$ is 
$(-1,-1)$, and around any $\alpha$ with $\alpha^{14}=1$ and $\alpha^2\not=1$
is $(-1,+1)$ (here $\bfz_2=\{+1,-1\}$). Hence the global monodromy group is 
the whole group $\bfz_2\times\bfz_2$. In particular $C$ is irreducible.\\

Notice also that the multiplicities of $f_2$ are all even numbers, but there 
is no germ $g:(X,x)\to (\bfc,0)$ with $f_2=g^2$.
(Actually, there is no homotopy between $f_2$ and $g^2$ for any $g$.)

\bekezdes{ Example.}{3.15}\ Set $(X,x)=(\{z^2+y(x^{12}-y^{18})=0\},0)$
and $f_1=x^2$ and $f_2=x^2-y^3$. Then $\G(X,f_1)= \G(X,f_2)$ is the graph:

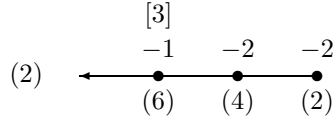
\begin{picture}(400,60)(0,0)
\put(110,20){\circle*{4}}
\put(140,20){\circle*{4}}
\put(170,20){\circle*{4}}
\put(110,20){\line(1,0){60}}
\put(110,30){\makebox(0,0){$-1$}}
\put(140,30){\makebox(0,0){$-2$}}
\put(170,30){\makebox(0,0){$-2$}}
\put(110,43){\makebox(0,0){$[3]$}}
\put(110,20){\vector(-1,0){30}}
\put(110,10){\makebox(0,0){$(6)$}}
\put(140,10){\makebox(0,0){$(4)$}}
\put(170,10){\makebox(0,0){$(2)$}}
\put(60,20){\makebox(0,0){$(2)$}}
\end{picture}

By a similar argument as in (\ref{3.14}) one has that 
  $arg_*(f_1)$  has cokernel $\bfz_2$, and $arg_*(f_2)$ is onto.

\bekezdes{ Example.}{3.16}\ 
Set $(X,x)=(\{z^2+(x^2-y^3)(x^3-y^2)=0\},0)$,
$f_1=x^2$ and $f_2=x^2-y^3$.
Then $\G(X,f_1)=\G(X,f_2)$ is the following graph:

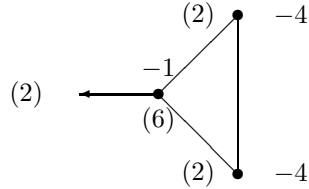
\begin{picture}(400,90)(0,-20)
\put(110,20){\circle*{4}}
\put(140,50){\circle*{4}}
\put(140,-10){\circle*{4}}
\put(110,20){\line(1,1){30}}
\put(110,20){\line(1,-1){30}}
\put(140,50){\line(0,-1){60}}
\put(110,30){\makebox(0,0){$-1$}}
\put(160,50){\makebox(0,0){$-4$}}
\put(125,50){\makebox(0,0){$(2)$}}
\put(160,-10){\makebox(0,0){$-4$}}
\put(125,-10){\makebox(0,0){$(2)$}}
\put(110,20){\vector(-1,0){30}}
\put(110,10){\makebox(0,0){$(6)$}}
\put(60,20){\makebox(0,0){$(2)$}}
\end{picture}

\noindent 
Then again:  $arg_*(f_1)$  has cokernel $\bfz_2$, and $arg_*(f_2)$ is onto.\\

Notice that in the above examples, $(X,x)=(\{z^2+h(x,y)=0\},0)$, and $f_2$ 
divides $h$ but it is not equal to $h$. 
For all such cases the monodromy argument given in (\ref{3.14}) is valid.
(So the interested reader can construct many--many  similar examples,
with even more additional properties.) But all these
 examples  define non--isolated singularities.
In order to construct examples of germs which define isolated singularities, 
we will use the well--known construction of series of singularities.
Namely, assume that $f_1$ and $f_2$ have the same graph but have different 
representations $arg_*$, and  their zero sets have non--isolated singularities 
(e.g. they are constructed by the above method).
Next, we find a germ $g$ such that the zero set of $f_i$ and $g$
have no common components (for $i=1,2$). 
Then, for a sufficiently large $k$
the germs $f_1+g^k$ and $f_2+g^k$
define isolated singularities whose embedded
resolution graphs are identical, but the representations $arg_*$ are 
different. 

\bekezdes{ Example.}{3.17}\  
Set $(X,x)=(\{x^2+y^7-z^{14}=0\},0)\subset
(\bfc^3,0)$ and take $f_1(x,y,z)=z^2$ and $f_2(x,y,z)=z^2-y$ as in (\ref{3.14}).
Let $P$ be the intersection point of the strict transform $S_a$ of
$\{f_i=0\}$ with the exceptional divisor $E$. Then, in some local  coordinate
system $(u,v)$ of $P$, $\{u=0\}$ represents $E$ (in a neighborhood of $P$), 
$\{v=0\}$ represents $S_a$, and $f_i=u^2v^2$. Consider $g=y$. Since $y$ in the
neighborhood of $P$ can be represented as $y=u^2$ (modulo a local invertible
germ), $f_i+g^k$ near $P$ has the form $u^2v^2+u^{2k}$. For example, if $k=2$,
then one needs one more blowing up in order to resolve $f_i+g^k$. 

Therefore, $\G(X,z^2+y^k)=\G(X,z^2-y+y^k)$ for any $k\geq 2$; and for 
$k=2$,  the graphs have the following form:

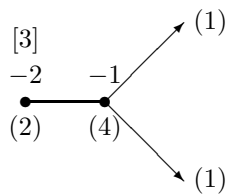
\begin{picture}(400,80)(0,-20)
\put(110,20){\circle*{4}}
\put(140,20){\circle*{4}}
\put(110,20){\line(1,0){30}}
\put(140,20){\vector(1,1){30}}
\put(140,20){\vector(1,-1){30}}

\put(110,30){\makebox(0,0){$-2$}}
\put(140,30){\makebox(0,0){$-1$}}
\put(110,43){\makebox(0,0){$[3]$}}
\put(110,20){\line(1,0){30}}

\put(110,10){\makebox(0,0){$(2)$}}
\put(140,10){\makebox(0,0){$(4)$}}
\put(180,50){\makebox(0,0){$(1)$}}
\put(180,-10){\makebox(0,0){$(1)$}}
\end{picture}

Notice that now ${\bf m}$ is onto, hence for both $i=1,2$, $arg_*(f_i)$
is onto. Nevertheless, $arg_*(f_1)\not=arg_*(f_2)$ because their restrictions
to a subgroup of $H_1(L_X\setminus L_f)$ are different. 

Indeed, 
let ${\caly}'$ be a tubular neighborhood of the irreducible exceptional
divisor $E$ of genus 3. This curve $E$ can be contracted by Grauert theorem
\cite{GRa}.  Then $E$ contracted in ${\caly}'$ gives birth to
a singularity $(X',x)$ with the same resolution graph as the surface 
singularity in (\ref{3.14}). Moreover, the germs $f_i$ ($i=1,2$) induce germs 
$f_i':(X',x)\to (\bfc,0)$, such that they 
have the same embedded resolution graphs as the germs in (\ref{3.14}). 
In particular,
$\coker (arg_*(f_1'))=\bfz_2$ and $arg_*(f_2')$ is onto. But $arg_*(f_i')$ is
the composed map 

\vspace{0.2cm}
\centerline{$H_1({\caly}'\setminus D,\bfz)\to H_1({\caly}\setminus D,\bfz)
\stackrel{arg_*(f_i)} {\longrightarrow}\bfz,$
}

\vspace{0.2cm}\noindent
hence $arg_*(f_1)\not=arg_*(f_2)$.

\bekezdes{ Example.}{3.18}\ 
Set $(X,x)=(\{z^2+(x^2-y^3)(x^3-y^2)=0\},0)$
and $f_1=x^2+y^k$ and $f_2=x^2-y^3+y^k$, where $k\geq 4$.
Then by similar argument as before, 
 $\G(X,f_1)=\G(X,f_2)$. This graph for $k=4$ is:

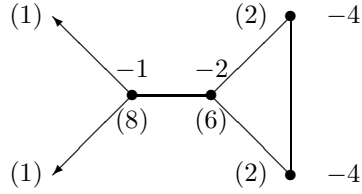
\begin{picture}(400,80)(-20,-20)
\put(80,20){\circle*{4}}
\put(110,20){\circle*{4}}
\put(140,50){\circle*{4}}
\put(140,-10){\circle*{4}}
\put(80,20){\line(1,0){30}}

\put(110,20){\line(1,1){30}}
\put(110,20){\line(1,-1){30}}
\put(140,50){\line(0,-1){60}}
\put(110,30){\makebox(0,0){$-2$}}
\put(80,30){\makebox(0,0){$-1$}}
\put(160,50){\makebox(0,0){$-4$}}
\put(125,50){\makebox(0,0){$(2)$}}
\put(160,-10){\makebox(0,0){$-4$}}
\put(125,-10){\makebox(0,0){$(2)$}}
\put(80,20){\vector(-1,-1){30}}
\put(80,20){\vector(-1,1){30}}
\put(110,10){\makebox(0,0){$(6)$}}
\put(80,10){\makebox(0,0){$(8)$}}
\put(40,50){\makebox(0,0){$(1)$}}
\put(40,-10){\makebox(0,0){$(1)$}}
\end{picture}

Let $E_i$ ($i=1,2,3$) be the irreducible exceptional divisors  with
self intersection numbers $-2,-4,-4$ respectively, and ${\caly}'$ be the 
union of small tubular neighborhoods of them. Then collapsing 
the curve $\cup_iE_i$ in
${\caly}'$ creates a singularity $(X',x)$ with the same graph as in (\ref{3.16}).
Repeating the arguments of (\ref{3.17}) (but using (\ref{3.16}) instead of 
(\ref{3.14}))
one has $arg_*(f_1)\not= arg_*(f_2)$. 

\subsection*{The universal cyclic covering of $\G(X,f)$.}

\vspace{2mm}

\bekezdes{-- The  covering $p:G(X,f)\to \G(X,f)$.}{5.1}

As we already noticed, the embedded resolution graph $\Gamma(X,f)$ does not
codify all the information about $arg_*$. On the other hand, this information
is needed in the study of cyclic coverings of $(X,x)$. 
In this section, we define a cyclic  covering of $\G(X,f)$ 
(cf. Section \ref{sec1}) 
which will control the behavior of the resolution graphs 
of all the cyclic coverings $\{X_{f,N}\}_N$
(cf. \ref{1.10}). The graph $G(X,f)$ was already considered in the literature
by Ph. Du Bois and F. Michel
from a completely different point of view, see \cite{MDB}. 

Let $(X,x)$ be  a normal surface singularity and $f:(X,x)\to (\bfc,0)$ the 
germ of an analytic function. Fix an embedded resolution $\phi:({\caly},D)
\to (X,f^{-1}(0))$ of $(f^{-1}(0),x)\subset (X,x)$ (as in \ref{1.1}) with  
embedded resolution graph $\G(X,f)$. 

Let $\tw$ ($w\in\calw$) be a small 
tubular neighborhood of the irreducible divisor $E_w$. 
By our assumption that any two irreducible exceptional divisor has at most one 
intersection point (see the first subsection of this chapter),
for any $e=(v,w)\in {\cale}$, 
the intersection $\tw\cap\tv$ ($(v,w)\in\calw\times\calw$) is
homeomorphic to a multidisc $D\times D$. This will be denoted by
$T_e$. If $T(S_a)$ ($a\in\cala$) is a small tubular neighborhood of 
the irreducible component $S_a$ of the strict transform $S$ (cf. \ref{1.1}),
 and
$a$ is adjacent to $w_a\in\calw$, then corresponding to the edge $e=(a,w_a)$
we introduce the multidisc $T_e=T(S_a)\cap T(E_{w_a})$. 
Set $T=(\cup_w\tw)\cup(\cup_aT(S_a))$.

Now, consider the smooth nearby fiber $f^{-1}(\delta)\subset X$ lifted  via 
$\phi$. For sufficiently small $\delta>0$, the fiber $F:=(f\circ\phi)^{-1}
(\delta)\subset {\caly}$ is in $T$. Set $F_w=F\cap \tw$ for any $w\in\calw$,
$F_a=F\cap T(S_a)$ for any $a\in\cala$, and $F_e=F\cap T_e$ for any $e\in\cale$.

It is not very difficult to construct a geometric
monodromy $h_g:F\to F$ of the 
fibration $f^{-1}(S_{\delta}^1)\to S^1_{\delta}$ (where $S^1_{\delta}=\{z
\in\bfc: |z|=\delta\}$ and $\delta$ is sufficiently small) which preserves the 
subspaces $\{F_v\}_{v\in\calv}$ and $\{F_e\}_{e\in\cale}$. Then the connected 
components of $F_v$ (resp. of $F_e$) are cyclically permuted by the geometric
monodromy. Let $n_v$ (resp. $n_e$) be the number of connected components
of $F_v$ (resp. $F_e$). Then, for any $e=(v_1,v_2)$, $n_e=d_e\cdot [n_{v_1},
n_{v_2}]$ for some $d_e\geq 1$. 

Now, we are able to construct the covering $p:G(X,f)\to \G(X,f)$
associated with the resolution $\phi$. 
Above a vertex
$v\in\calv(\G(X,f))$ there are exactly $n_v$ vertices of $G(X,f)$, they 
correspond to the connected components of $F_v$. The $\bfz$--action is induced 
by the monodromy (by the identification $1_{\bfz}=(h_g)_*$).
If $v$ is an arrowhead in $\G$ then by convention, all the
 vertices in $G$ above $v$
are arrowheads (cf. \ref{4.20} (1)). 
Above an edge $e$ of $\G$, there are $n_e$ edges of 
$G$. They corresponds to the connected components of $F_e$. The $\bfz$--action
is again  generated by the monodromy. 

Fix an edge $\tilde{e}$ of $G$ (above the edge $e$ of $\G$) which corresponds
 to the connected component $F_{\tilde{e}}$ of $F_e$. Similarly, take a vertex
$\tilde{v}$ of $G$ (above the vertex $v$ of $\G$) which corresponds to the 
connected component $F_{\tilde{v}}$ of $F_v$. Then $\tilde{e}$ has as an
 endpoint  the vertex $\tilde{v}$ (in $G$) if and only if $F_{\tilde{e}}
\subset F_{\tilde{v}}$. In particular, with the same notations as above,
$\tilde{v_1}$ and $\tilde{v_2}$ are connected in $G$ by (at least) one
edge if and only if $F_{\tilde{v}_1}\cap F_{\tilde{v}_2}\not=\emptyset$. 
Set $e=(v_1,v_2)$. 
If $F_{\tilde{v}_1}\cap F_{\tilde{v}_2}\not=\emptyset$, then 
$F_{\tilde{v}_1}\cap F_{\tilde{v}_2}$ has exactly $d_e=[n_{v_1},
n_{v_2}]/n_e$  connected components, hence $\tilde{v}_1$ and 
$\tilde{v}_2$ are connected exactly by $d_e$ edges. Therefore, 
above the segments of $\G$ we have exactly the  ``standard
blocks'' of (\ref{4.2}) in $G$.

In the sequel, we will use the notation $({\bf n},{\bf d})=
\{\{n_v\}_{v\in\calv};\{d_e\}_{e\in\cale}\}$, where $n_e=d_e[n_{v_1},n_{v_2}]$
for any edge $e=(v_1,v_2)$. 
 
\vspace{2mm}

The next lemma establishes the number of connected components of $G(X,f)$.

\bekezdes{Lemma.}{5.3} {\em The number of connected components of the 
graph $G(X,f)$ is equal to the number of connected components of the Milnor
fiber $F$ of the germ $f$.}

\noindent 
{\em Proof.}\ Let $|G|$ be the topological realization of the graph $G$
considered as a 1--dimensional simplicial complex. Then it is not difficult 
to construct a continuous map $\alpha:F\to |G|$ which maps $F_{\tilde{v}}$
 to the zero-cell (vertex) $\tilde{v}$, maps $F_{\tilde{e}}$ to the one--cell
(edge) $\tilde{e}$, and for any $P\in|G|$, the space $\alpha^{-1}(P)$ is 
connected. (Notice that $F_{\tilde{e}}\approx S^1\times \tilde{e}$, then
$\alpha$ restricted to $F_{\tilde{e}}$ can be identified with the second 
projection $S^1\times \tilde{e}\to\tilde{e}$.) Hence $\pi_0(\alpha):
\pi_0(F)\to \pi_0(|G|)$ is an isomorphism.\hfill $\diamondsuit$

\vspace{2mm}

Now, it is well--known, that the fibrations $f:f^{-1}(S^1_{\delta})\to 
S^1_{\delta}$ and $arg=f/|f|:L_X\setminus L_f\to S^1$ are equivalent.
Hence by the  long homotopy exact sequence:
$$\pi_1(L_X\setminus L_f)\stackrel{\pi_1(arg)}{\to}\pi_1(S^1)=\bfz\to \pi_0(F)
\to 0$$
(where $\pi_1(arg)=ab\circ arg_*$,
cf. \ref{3.10}) we obtain that $|\pi_0(F)|=|\coker(arg_*)|$. Therefore:

\bekezdes{Corollary.}{5.4} {\em 
 $G(X,f)$ has $|\coker(arg_*(f))|$ connected components. In particular,
if $f$ defines an isolated singularity, then $G(X,f)$ is a connected graph.
(But $G(X,f)$ may be connected  even for germs $f$ with
$g.c.d.\{m_v:v\in\calv\}\not=1$.)}

\vspace{0.2cm}\noindent
Before we state the second part of this corollary, we make the following 
discussion. 

Fix  a connected  subgraph $\G'$  of $\G(X,f)$
with non--arrowhead vertices $\calw'$.  Since the intersection form
associated with the exceptional divisors 
$E':=\cup_{w\in\calw'}E_w$ is negative definite, by Grauert theorem \cite{GRa}
$E'\subset {\caly}$ 
can be contracted to a singular point. Let $(X',x')$ be this singular point.
Moreover, since $f\circ\phi$ is 
zero along $E'$, it gives rise to a germ $f'$ defined on $(X',x')$.
It is obvious, that $\G'$ is the  embedded 
resolution graph of a singularity $((f')^{-1}(0),x')\subset (X',x')$.
 Let the 
corresponding representation be denoted by $arg_*(f')$. 

\bekezdes{Corollary.}{5.4b} {\em 

Consider the universal cyclic covering  $p:G(X,f)\to
\G(X,f)$ associated with $f$ and the resolution $\phi$.
Fix a connected  subgraph $\G'$  of $\G(X,f)$. Then
$p^{-1}(\G')$ has $|\coker(arg_*(f'))|$ connected components.}

\bekezdes{-- The case when $L_X$ is a rational homology sphere.}{5.5}
If $L_X$ is a rational homology sphere, then the Milnor fibration is
completely determined by ${\bf m}$ (cf. \ref{3.13}). Hence,
$G(X,f)$ contains the same amount of information as $\G(X,f)$, and it can 
always be recovered from $\G(X,f)$. 

\bekezdes{Lemma.}{L}\  {\em Assume that $L_X$ is a rational homology sphere.
Then $n_v:= g.c.d.\{m_w:\ w\in\calv_v\cup\{v\}\}$
for any $v\in\calv(\G)$; 
and $n_e:=g.c.d.(m_{v_1},m_{v_2})$
for any $e=(v_1,v_2)\in\cale(\G)$.
(In particular, for $a\in\cala(\G)$ one has
$n_a=n_e$, where $e=(a,w_a)\in\cale(\G)$.) Moreover, the number of connected 
components of $G(X,f)$ is exactly $g.c.d.\{m_v: v\in\calv(\G)\}$.}

\vspace{2mm}

\noindent {\em Proof.}\ If $T_v$ is a tubular neighborhood of $E_v$ as in 
(\ref{5.1}), then $\coker(\pi_1(T_v\setminus \phi^{-1}(f^{-1}(0)))\to\bfz)$
is $\bfz/\{m_w:w\in\calv_v\cup\{v\}\}\bfz$, hence the statement about $n_v$
follows (by \ref{5.4b}). Similarly, for $e=(v,w)$, the cokernel of 
$\pi_1(T_v\cap T_w\setminus \phi^{-1}(f^{-1}(0)))\to\bfz)$ is $\bfz/(m_v\bfz+
m_w\bfz)$. The number of connected components of $G$ is 
 $|\coker(arg_*)|=|\coker ({\bf m}:H_{\G}\to\bfz)|$. \hfill $\diamondsuit$

\vspace{2mm}

Notice also that in this case  ${\calg}(\G(X,f),({\bf n},{\bf d}))$ 
contains only one element
(because $\G$ is a tree, cf. \ref{4.16}). This class is represented by $G(X,f)$.\\

For example, if $(X,x)$ is smooth, then $L_X=S^3$, hence for 
any plane curve singularity, the universal cyclic
covering $G(X,f)\to \G(X,f)$ can  completely be determined from 
the embedded resolution graph $\G(X,f)$ of $f$.

\bekezdes{Corollary.}{5.7} {\em 
Consider an arbitrary germ $f:(X,x)\to (\bfc,0)$ (i.e. without any 
restriction about $L_X$), and consider the universal cyclic
covering $p:G(X,f)\to
\G(X,f)$. Then for any $w\in\calw(\G)$ with $g_w=0$, $p^{-1}(w)$ 
consists of exactly $n_w=g.c.d.\{m_v:v\in\calv_w\cup\{w\}\}$  vertices
of $G$.
Similarly, for any $a\in\cala(\G)$, $\#p^{-1}(a)=n_a=g.c.d.(m_a,m_{w_a})$,
where $(a,w_a)\in\cale(\G)$. Moreover, the number of edges $n_e$ above
$e=(v_1,v_2)$ is $g.c.d.(m_{v_1},m_{v_2})$. 

For $w\in\calw(\G)$ with $g_w>0$ the following divisibility holds:
\ $n_w\,|\,g.c.d.\{m_v:v\in\calv_w \cup\{w\}\}$.}

\vspace{2mm}

\noindent 
{\em Proof.}\ Notice that if an irreducible rational exceptional divisor 
is contracted to a singular point then the link of this singular point is
a rational homology sphere. Then use (\ref{5.4}--\ref{L}).\hfill 
$\diamondsuit$

\bekezdes{-- The case when $\G(X,f)$ is a tree.}{5.8} 
Recall that if $\G(X,f)$ is a tree, then for any system of integers 
$({\bf n},{\bf d})$,
 by (\ref{4.16}) the class ${\calg}(\G(X,f),({\bf n},{\bf d}))=0$, 
hence all the coverings are equivalent.
One the other hand, if $L_X$ is not rational homology sphere, even if 
$\G(X,f)$ is a tree, the covering data  $\{n_x\}_{x\in\calv\cup\cale}$ of the 
universal covering $G\to \G$ (more precisely, the integers $n_w$ with 
$g_w>0$) are not determined by $\G(X,f)$ (see the next examples).
Hence already in this case, the covering $G(X,f)\to \G(X,f)$ contains 
some additional information about the representation $arg_*$.

\bekezdes{Example.}{5.9} 
Set $(X,x)=(\{x^2+y^7-z^{14}=0\},0)\subset
(\bfc^3,0)$ and take $f_1(x,y,z)=z^2$ and $f_2(x,y,z)=z^2-y$ (cf. \ref{3.14}). 
Then the coverings $p:G(X,f_i)\to \G(X,f_i)$ (for $i=1,2$) are:
 
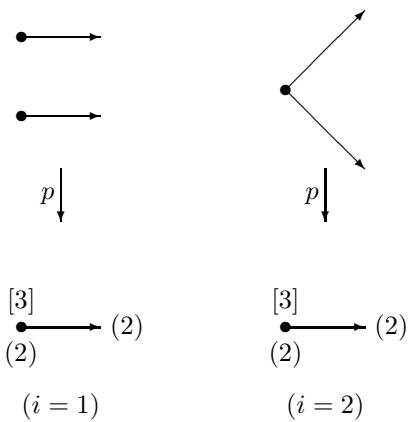
\hspace{-1cm}
\begin{picture}(400,180)(0,-20)
\put(110,20){\circle*{4}}
\put(110,30){\makebox(0,0){$[3]$}}
\put(110,20){\vector(1,0){30}}
\put(110,10){\makebox(0,0){$(2)$}}
\put(150,20){\makebox(0,0){$(2)$}}
\put(125,-10){\makebox(0,0){$(i=1)$}}
\put(125,80){\vector(0,-1){20}}
\put(120,70){\makebox(0,0){$p$}}

\put(110,100){\circle*{4}}
\put(110,130){\circle*{4}}
\put(110,100){\vector(1,0){30}}
\put(110,130){\vector(1,0){30}}

\put(210,20){\circle*{4}}
\put(210,30){\makebox(0,0){$[3]$}}
\put(210,20){\vector(1,0){30}}
\put(210,10){\makebox(0,0){$(2)$}}
\put(250,20){\makebox(0,0){$(2)$}}
\put(225,-10){\makebox(0,0){$(i=2)$}}
\put(225,80){\vector(0,-1){20}}
\put(220,70){\makebox(0,0){$p$}}

\put(210,110){\circle*{4}}
\put(210,110){\vector(1,1){30}}
\put(210,110){\vector(1,-1){30}}
\end{picture}

\noindent In order to count the number of vertices above the irreducible 
divisor $E$  (with $g=3$), we have to consider the representations $arg_*(f_i)$
(cf. \ref{5.4}). In the first case, the cokernel of this representation is
${\bfz}_2$, in the second case it is trivial. Hence, above $E$, in the first
case one has two vertices, and in the second case only one.
 
\vspace{0.1cm} \noindent
Notice, that the number of connected components (and even the 
Euler--character-\\
istic) of the graphs $G(X,f_i)$ are different.

\bekezdes{ Example.}{5.10} 
Set $(X,x)=(\{x^2+y^7-z^{14}=0\},0)\subset
(\bfc^3,0)$ and take $f_1(x,y,z)=z^2+y^2$ and $f_2(x,y,z)=z^2-y+y^2$ 
(cf. \ref{3.17}).
By a similar argument as above, 
the coverings $p:G(X,f_i)\to \G(X,f_i)$ (for $i=1,2$) are:

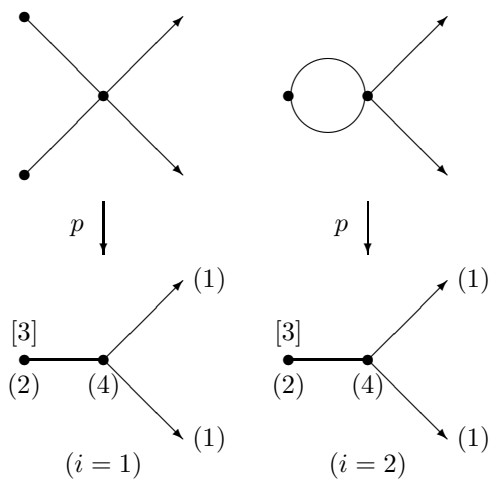
\hspace{-1cm}
\begin{picture}(400,200)(0,-40)
\put(110,20){\circle*{4}}
\put(140,20){\circle*{4}}
\put(110,20){\line(1,0){30}}
\put(140,20){\vector(1,1){30}}
\put(140,20){\vector(1,-1){30}}
\put(110,30){\makebox(0,0){$[3]$}}
\put(110,10){\makebox(0,0){$(2)$}}
\put(140,10){\makebox(0,0){$(4)$}}
\put(180,50){\makebox(0,0){$(1)$}}
\put(180,-10){\makebox(0,0){$(1)$}}
\put(140,-20){\makebox(0,0){$(i=1)$}}

\put(140,80){\vector(0,-1){20}}
\put(130,70){\makebox(0,0){$p$}}

\put(110,90){\circle*{4}}
\put(110,150){\circle*{4}}
\put(140,120){\circle*{4}}
\put(110,90){\line(1,1){30}}
\put(110,150){\line(1,-1){30}}
\put(140,120){\vector(1,1){30}}
\put(140,120){\vector(1,-1){30}}

\put(210,20){\circle*{4}}
\put(240,20){\circle*{4}}
\put(210,20){\line(1,0){30}}
\put(240,20){\vector(1,1){30}}
\put(240,20){\vector(1,-1){30}}
\put(210,30){\makebox(0,0){$[3]$}}
\put(210,10){\makebox(0,0){$(2)$}}
\put(240,10){\makebox(0,0){$(4)$}}
\put(280,50){\makebox(0,0){$(1)$}}
\put(280,-10){\makebox(0,0){$(1)$}}
\put(240,-20){\makebox(0,0){$(i=2)$}}

\put(240,80){\vector(0,-1){20}}
\put(230,70){\makebox(0,0){$p$}}

\put(210,120){\circle*{4}}
\put(240,120){\circle*{4}}
\put(225,120){\circle{29}}
\put(240,120){\vector(1,1){30}}
\put(240,120){\vector(1,-1){30}}

\end{picture}

\noindent In this case the number of independent cycles 
of the graphs $G(X,f_i)$ is different.

\bekezdes{}{5.11} 
If all the irreducible exceptional divisors of $\phi$  are rational (i.e. 
$g_w=0$ for all $w$), then the type  $({\bf n},{\bf d})$ of the covering
 is completely
determined by $\G(X,f)$ (cf. \ref{5.7}).
But if $\G(X,f)$ is not a tree, then ${\calg}(\G(X,f),({\bf n},{\bf d}))$
can be non--trivial. Hence again, $G(X,f)$ carries some additional 
information about  $arg_*$. 

\bekezdes{ Example.}{5.12} 
Set $(X,x)=(\{z^2+(x^2-y^3)(x^3-y^2)=0\},0)$
and $f_1=x^2+y^4$ and $f_2=x^2-y^3+y^4$ (cf. \ref{3.18}). 
Then the coverings $p:G(X,f_i)\to \G(X,f_i)$ (for $i=1,2$) are:

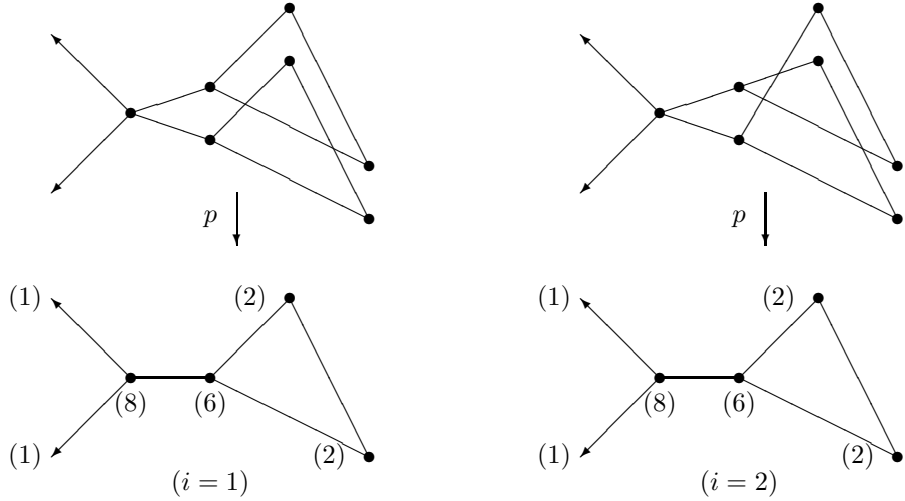
\hspace{-1cm}
\begin{picture}(400,210)(20,-30)
\put(80,20){\circle*{4}}
\put(110,20){\circle*{4}}
\put(140,50){\circle*{4}}
\put(170,-10){\circle*{4}}
\put(80,20){\line(1,0){30}}
\put(110,-20){\makebox(0,0){$(i=1)$}}
\put(110,20){\line(1,1){30}}
\put(110,20){\line(2,-1){60}}
\put(140,50){\line(1,-2){30}}
\put(125,50){\makebox(0,0){$(2)$}}
\put(155,-10){\makebox(0,0){$(2)$}}
\put(80,20){\vector(-1,-1){30}}
\put(80,20){\vector(-1,1){30}}
\put(110,10){\makebox(0,0){$(6)$}}
\put(80,10){\makebox(0,0){$(8)$}}
\put(40,50){\makebox(0,0){$(1)$}}
\put(40,-10){\makebox(0,0){$(1)$}}

\put(120,90){\vector(0,-1){20}}
\put(110,80){\makebox(0,0){$p$}}

\put(80,120){\circle*{4}}
\put(110,110){\circle*{4}}
\put(110,130){\circle*{4}}
\put(140,160){\circle*{4}}
\put(140,140){\circle*{4}}
\put(170,100){\circle*{4}}
\put(170,80){\circle*{4}}
\put(80,120){\line(3,1){30}}
\put(80,120){\line(3,-1){30}}
\put(110,110){\line(1,1){30}}
\put(110,130){\line(1,1){30}}
\put(110,110){\line(2,-1){60}}
\put(110,130){\line(2,-1){60}}
\put(140,160){\line(1,-2){30}}
\put(140,140){\line(1,-2){30}}
\put(80,120){\vector(-1,-1){30}}
\put(80,120){\vector(-1,1){30}}

\put(280,20){\circle*{4}}
\put(310,20){\circle*{4}}
\put(340,50){\circle*{4}}
\put(370,-10){\circle*{4}}
\put(280,20){\line(1,0){30}}
\put(310,-20){\makebox(0,0){$(i=2)$}}
\put(310,20){\line(1,1){30}}
\put(310,20){\line(2,-1){60}}
\put(340,50){\line(1,-2){30}}
\put(325,50){\makebox(0,0){$(2)$}}
\put(355,-10){\makebox(0,0){$(2)$}}
\put(280,20){\vector(-1,-1){30}}
\put(280,20){\vector(-1,1){30}}
\put(310,10){\makebox(0,0){$(6)$}}
\put(280,10){\makebox(0,0){$(8)$}}
\put(240,50){\makebox(0,0){$(1)$}}
\put(240,-10){\makebox(0,0){$(1)$}}

\put(320,90){\vector(0,-1){20}}
\put(310,80){\makebox(0,0){$p$}}

\put(280,120){\circle*{4}}
\put(310,110){\circle*{4}}
\put(310,130){\circle*{4}}
\put(340,160){\circle*{4}}
\put(340,140){\circle*{4}}
\put(370,100){\circle*{4}}
\put(370,80){\circle*{4}}
\put(280,120){\line(3,1){30}}
\put(280,120){\line(3,-1){30}}
\put(310,110){\line(3,5){30}}

\put(310,130){\line(3,1){30}}
\put(310,110){\line(2,-1){60}}
\put(310,130){\line(2,-1){60}}
\put(340,160){\line(1,-2){30}}
\put(340,140){\line(1,-2){30}}
\put(280,120){\vector(-1,-1){30}}
\put(280,120){\vector(-1,1){30}}

\end{picture}

\noindent
In this case, ${\calg}(\G(X,f), ({\bf n},{\bf d}))=\bfz_2$ (cf. \ref{4.18}). The above 
examples provide the two non--equivalent representatives of ${\calg}=\bfz_2$.

The number of independent cycles $c_G$ for both graphs $G(X,f_i)$ are the same
($=2$). This is a general fact: for a connected graph $G$ one has: 
$1-c_G=\#\calv(G)-\#\cale(G)$, but $\#\calv(G(X,f_i))$ and $\#\cale(G(X,f_i))$
are determined  by $\G$ (provided that $g_w=0$ for any $w$).
 
\section{The resolution graph of cyclic coverings}\label{sec3}

\subsection*{The monodromy representation  of cyclic coverings.}

\vspace{2mm}

\bekezdes{}{6.1} \ Let $(X,x)$ be a normal surface singularity and 
$f: (X,x)\to (\bfc,0)$ the germ of an analytic function.
For any integer $N\geq 1$, take $b: (\bfc , 0)\to (\bfc ,0)$
given by $z\longmapsto z^{N}$, and let $X_{f,N}$ be the normalization
of the fiber product $\{ (x',z)\in (X\times\bfc , x\times 0) : f(x')=z^{N}\}$
(cf. \ref{1.10}).
The second projection $(x',z)\in X\times\bfc\longmapsto z\in\bfc$ induces an analytic map 
$X_{f,N}\to\bfc$, still denoted by $z$.
The first projection $(x',z)\longmapsto x'$ gives 
rise to a ramified cyclic $N$-covering 
$pr  :X_{f,N}\to X$, branched along $\finv$. If $\finv$ has an isolated
 singular point at $x$, 
then there is only one (singular) point of $\xfn$ lying above $x\in X$.
But in general $pr ^{-1}(x)$ contains more than one point (cf. \ref{1.7}). 

We regard  $\bfz _N$ as the group of $N^{th}$-roots of unity 
$\{ \xi_k=e^{2\pi ik/N};$ $\ 0\le k\le N-1\}$,
then $(x',z)\longmapsto (x',\xi_k z)$ induces a 
$\bfz _N$--Galois action of $\xfn$ over $X$.

If $P\in X\setminus\finv$ is a point in the complement of the branch locus, then
$pr ^{-1}(P)$ consists of $N$ points and they are cyclically permuted by the Galois
action of $\bfz _N$. The monodromy representation of the regular covering
$pr  \vert _{\xfn\setminus \{ z=0 \}} : \xfn\setminus \{ z=0\}\to X\setminus\{ f=0\}$
is denoted by $\varphi_N:\pi_1(X\setminus \{ f=0\})\to\bfz_N$.
Notice that $X\setminus \{f=0\}$ is connected, and $\bfz_N$ is abelian, hence
the monodromy representation does not depend on the choice of the basepoint.
So, we will omit the basepoint of the fundamental group.

By the local cone structure of $(X,\{ f=0\} )$, we can replace the group $\pi_1 (X\setminus\{ f=0\} )$
by $\pi_1 ( L_X\setminus L_f)$. 
The following property of  the cyclic coverings
is well-known:

\bekezdes{Lemma.}{6.2}\  {\em 
Let $\pi_1(  arg ): \pi_1 (L_X\setminus L_f)\to\pi_1(S^1)=\bfz$ 
be the morphism induced by the Milnor fibration of $f$ (at the fundamental 
group level) (cf. \ref{3.10}); and let
$mod_N:\bfz\to\bfz_N$ be the natural projection $1\longmapsto \hat{1}$.
Then $\varphi_N=mod_N\circ\pi_1 (  arg )$.}

\vspace{2mm}

Since $\bfz_N$ is abelian, $\varphi_N = r_N\circ ab$ for some 
$r_N: H_1(L_X\setminus L_f,\bfz )\to\bfz_N$ (where $ab :\pi_1\to
\pi_1/[\pi_1,\pi_1]=
H_1$ is the abelianization map).
Hence, $r_N =mod_N\circ arg_*$, where $ arg_*: H_1(L_X\setminus L_f,\bfz )\to\bfz $ is 
induced by the Milnor fibration $ arg = f/|f|$. This implies the following:

\bekezdes{Lemma.}{6.3} 
{\em If the Milnor fiber $F$ of $f$ has $k$ connected components, then
$\xfn$ has $(k,N)$ connected components.}

\vspace{2mm}

\noindent
{\em Proof.}\ By the 
long homotopy exact sequence of the Milnor fibration $k=|\coker arg_*|$,
 hence the range of  $arg_*$ is $k\bfz\subset\bfz$.
But $mod_N(k\bfz )=k\bfz_N$ has index $(k,N)$ in $\bfz_N$, hence 
$\coker r_N=\coker\varphi_N$ $\approx\bfz_{(k,N)}$. Again, by the long homotopy exact sequence of the
regular covering $pr $ over $X\setminus \{ f=0\}$, 
the integer $|\coker\varphi_N|$ is the number of connected 
components of $\xfn\setminus \{ z=0\}$.
But this is exactly the number of connected (or irreducible) components of $\xfn$.
\hfill$\diamondsuit$

\vspace{2mm}

\noi
The number $(N,k)$ is exactly the number of points of $\xfn$ lying above $x\in X$ 
(i.e. $\#pr ^{-1}(x)$). In the sequel, we will use the germ--notation 
$(\xfn ,\{x_1,...,x_{(N,k)}\})$, which means that $\xfn$ consists of $(N,k)$ 
disjoint space germs
$(\xfn , x_i)_{i=1}^{(N,k)}$. Obviously, they are all isomorphic with each other -- an
isomorphism is given by the Galois action (which permutes the points $\{ x_i\}_{i=1}^{(N,k)}$).

\bekezdes{Remark.}{6.4}
The number of (singular) points of $\xfn$ lying above $x\in X$ cannot be determined from the
embedded resolution graph of $f$ and the integer $N$. 

Indeed, consider the situation described in (\ref{3.14}) 
and (\ref{5.9}). Then 
$\#pr ^{-1}(x)=2$
in the first case, and  $\#pr ^{-1}(x)=1$ in the second case.

\vspace{2mm}

The above remark already suggests (and we will see a lot of other examples
later) that  the resolution graph of $\xfn$
cannot be reconstructed from the graph $\Gamma (X,f)$ and $N$.
 In fact, this was the very reason why
 we constructed the universal cyclic covering $G(X,f)\to\Gamma (X,f)$.
For example, related to the above discussion: 
the number of connected 
components of $G(X,f)$ is exactly $k$ (cf. \ref{5.3}), hence 
$(k,N)=\#pr ^{-1}(x)$ is determined by $G(X,f)$ and the integer $N$. 

\bekezdes{Definitions.}{6.5}\ 

{\em
a.) The resolution graph $\Gamma (\xfn )$, by definition, is the union of the resolution graphs
of $(\xfn ,x_i)_{i=1}^{(k,N)}$.

b.) The embedded resolution graph $\Gamma (\xfn ,z)$ of
$z:(\xfn ,\{x_1,...,x_{(k,N)}\} )$ $\to (\bfc,0)$ is the union of the embedded resolution 
graphs of $z:(\xfn ,x_i )\to(\bfc,0)\ \ (1\le i\le (k,N))$.}

\vspace{2mm}

In both cases $\G$ has $(k,N)$ identical connected components.

\subsection*{The embedded resolution graph of  $\G(X_{f,N},z)$.}

\vspace{2mm}

\bekezdes{ -- The main construction.}{7.1}  Let $p:G(X,f)\to \G(X,f)$
be the universal cyclic covering  of the embedded resolution graph
$\G(X,f)$  associated with the germ $f:(X,x)\to (\bfc,0)$ (cf. \ref{5.1}).
For brevity we will use $\G=\G(X,f)$ and $G=G(X,f)$ throughout this 
subsection. 
The covering $G$ is an element of $\fedog$, where $n_v=\#p^{-1}(v)\subset 
\calv(G)$ for any vertex $v\in\cvg$, and $n_e=d_e[n_{v_1},n_{v_2}]=\#
p^{-1}(e)\subset \cale(G)$ for any $e\in\evg$.
Recall that the graph $\G$ has the following  decorations: multiplicities
$\{m_v\}_{v\in\cvg}$, genera $[g_w]_{w\in\calw(\G)}$ 
(and self intersection numbers, which are less important in this construction).
By our convention: $\cala(G)=p^{-1}(\cala(\G))$.

Now, for any fixed integer $N\geq 1$,
we construct a new graph in four steps.\\

\noindent 
{\bf Step 1.}\ The graph $G$ has a $\bfz$--action. The ``orbit
graph'' of the subgroup $N\bfz\subset \bfz$ is denoted by $mod_N(G)$
(for details, see \ref{4.19}). \\
The new covering $mod_N(p):
mod_N(G)\to\G$ is an element of
${\calg}(\G,((N,{\bf n}),{\bf d'}))$, where $(N,{\bf n})=
\{g.c.d.(N,n_v)\}_{v\in\cvg}$, and $d_e'=(d_e,N/(N,[n_{v_1},n_{v_2}]))$ for any
edge $e=(v_1,v_2)\in\evg$ (i.e. $\#mod_N(p)^{-1}(e)=g.c.d.(N,\#p^{-1}(e))$
for any $e$).\\

\noindent {\bf Step 2.}\ We put decorations 
(multiplicities and genera) on $mod_N(G)$ as follows. 

(a)\ The multiplicity $(m_{\tilde{v}})$ of any vertex $\tilde{v}\in\calv(
mod_N(G))$, which lies above $v\in\cvg$ is $m_v/(m_v,N)$;

(b)\ The genus $[\tilde{g}_{\tilde{w}}]$ of any vertex $\tilde{w}\in
\calw(mod_N(G))$, which lies above $w\in\calw(\G)$ with genus $g_w$, is
given by:
$$2-2\tilde{g}_{\tilde{w}}=\frac{(2-2g_w-\delta_w)\cdot (m_w,N)+\sum_{v\in\calv_w(\G)}g.c.d.(m_v,m_w,N)}{(N,n_w)},$$
where $\delta_w=\#\calv_w(\G)$.\\

\noindent {\bf Step 3.}\ Any edge $\tilde{e}$ of $mod_N(G)$, with endpoints
$\tilde{v_1}$ and $\tilde{v_2}$, 
lying above $e$ with endpoints $v_1$ and $v_2\in\evg$,  
will be replaced by a string 
$\overline{Str}(e)$ as follows (cf. \ref{2.3} and \ref{4.20} (3)). 

If $t_e=g.c.d.(m_{v_1},m_{v_2},N)$, then set $N'=n/t_e$ and $m_{v_i}'=m_{v_i}/
t_e$ for $i=1,2$. 
Consider  the unique $0\leq \lambda <N'/(m_{v_1}',N')$  and 
$m_1\in{\bf N}$ with:
$$m_{v_2}'+\lambda \cdot \frac{m_{v_1}'}{(m_{v_1}',N')}=m_1\cdot 
\frac{N'}{(m_{v_1}',N')}.$$
If $\lambda=0$, then the edge $\tilde{e}$ remains unchanged.\\
If $\lambda\not=0$, then take  the continuous fraction:
$$\frac{N'/(m_{v_1}',N')}{\lambda}=k_1-{1\over\displaystyle
k_2-{\strut 1\over\displaystyle\ddots
-{\strut 1\over k_s}}}, \ \ k_1,\ldots, k_s\geq 2.$$
If both $v_1$ and $v_2$ are non--arrowhead vertices, then 
$\overline{Str}(e)$ denotes the following
decorated string:

\begin{picture}(400,80)(30,0)
\put(150,35){\makebox(0,0){$-k_1$}}
\put(200,35){\makebox(0,0){$-k_2$}}
\put(300,35){\makebox(0,0){$-k_s$}}
\put(150,15){\makebox(0,0){$(m_1)$}}
\put(100,15){\makebox(0,0){$(m_{\tilde{v}_1})$}}
\put(350,15){\makebox(0,0){$(m_{\tilde{v}_2})$}}
\put(200,15){\makebox(0,0){$(m_2)$}}
\put(300,15){\makebox(0,0){$(m_s)$}}
\put(150,50){\makebox(0,0){$[0]$}}
\put(100,50){\makebox(0,0){$[\tilde{g}_{\tilde{v}_1}]$}}
\put(350,50){\makebox(0,0){$[\tilde{g}_{\tilde{v}_2}]$}}
\put(200,50){\makebox(0,0){$[0]$}}
\put(300,50){\makebox(0,0){$[0]$}}
\put(30,25){\makebox(0,0){$\overline{Str}(e):$}}

\put(150,25){\circle*{4}}
\put(100,25){\circle*{4}}
\put(350,25){\circle*{4}}
\put(200,25){\circle*{4}}
\put(300,25){\circle*{4}}
\put(225,25){\line(-1,0){125}}
\put(275,25){\line(1,0){75}}
\put(250,25){\makebox(0,0){$\cdots$}}
\end{picture}

\noindent 
with  genera $[\tilde{g}_{\tilde{v}_1}],\ [0],\ \ldots, [0],\ [\tilde{g}_{\tilde{v}_2}]$, self intersection numbers $-k_1,\ldots,-k_s$, and multiplicities
$m_{\tilde{v}_1},\ m_1,\ \ldots, m_s,\ m_{\tilde{v}_2}$. 

The multiplicities $m_{\tilde{v}_1}$ and $m_{\tilde{v}_2}$
were already determined in step 2,
namely $m_{\tilde{v}_i}=m_{v_i}/(m_{v_i},N)$; 
and  $m_1$ is the number given by the above congruence. 
Moreover,  the multiplicities $m_2,\ldots, m_s$
can  easily be determined using (\ref{1.3}\  (2)), namely
$$m_2=k_1m_1-m_{\tilde{v}_1}; \ \mbox{and}\
m_{i+1}=k_im_i-m_{i-1}\ \mbox{for $i\geq 2$}.$$

\noindent 
Then each edge $\tilde{e}$:

\begin{picture}(400,50)(30,0)
\put(100,15){\makebox(0,0){$(m_{\tilde{v}_1})$}}
\put(150,15){\makebox(0,0){$(m_{\tilde{v}_2})$}}
\put(100,35){\makebox(0,0){$[\tilde{g}_{\tilde{v}_1}]$}}
\put(150,35){\makebox(0,0){$[\tilde{g}_{\tilde{v}_2}]$}}
\put(150,25){\circle*{4}}
\put(100,25){\circle*{4}}
\put(100,25){\line(1,0){50}}
\end{picture}

\noindent of $mod_N(G)$, lying above $e$, 
is replaced by the string $\overline{Str}(e)$. 

Moreover, if $v_2\in\cala(\G)$, then
the  edge $\tilde{e}$:

\begin{picture}(400,50)(30,0)
\put(100,15){\makebox(0,0){$(m_{\tilde{v}_1})$}}
\put(100,35){\makebox(0,0){$[\tilde{g}_{\tilde{v}_1}]$}}
\put(180,25){\makebox(0,0){$(m_{\tilde{v}_2})$}}
\put(100,25){\circle*{4}}
\put(100,25){\vector(1,0){50}}
\end{picture}

\noindent of $mod_N(G)$, lying above $e$, 
is replaced by the ``modified'' string $\overline{Str}(e)$
(which has the same decorations as the ``original $\overline{Str}(e)$''):

\hspace{-1.5cm}
\begin{picture}(400,80)(30,0)
\put(150,35){\makebox(0,0){$-k_1$}}
\put(200,35){\makebox(0,0){$-k_2$}}
\put(300,35){\makebox(0,0){$-k_s$}}
\put(150,15){\makebox(0,0){$(m_1)$}}
\put(100,15){\makebox(0,0){$(m_{\tilde{v}_1})$}}
\put(380,25){\makebox(0,0){$(m_{\tilde{v}_2})$}}
\put(200,15){\makebox(0,0){$(m_2)$}}
\put(300,15){\makebox(0,0){$(m_s)$}}
\put(150,50){\makebox(0,0){$[0]$}}
\put(100,50){\makebox(0,0){$[\tilde{g}_{\tilde{v}_1}]$}}
\put(200,50){\makebox(0,0){$[0]$}}
\put(300,50){\makebox(0,0){$[0]$}}

\put(150,25){\circle*{4}}
\put(100,25){\circle*{4}}
\put(200,25){\circle*{4}}
\put(300,25){\circle*{4}}
\put(225,25){\line(-1,0){125}}
\put(275,25){\vector(1,0){75}}
\put(250,25){\makebox(0,0){$\cdots$}}
\end{picture}

\noindent The new graph resulting from inserting all
the necessary strings into $mod_N(G)$ is denoted by $mod_N(G)({\bf Str})$.\\

\noindent {\bf Step 4.}\ The decoration of $mod_N(G)({\bf Str})$ is not 
complete. All the vertices have multiplicities, all the non--arrowhead
vertices have genera, but some of the self intersection numbers are
missing (corresponding exactly to the vertices $\tilde{v}$ of 
$mod_N(G)$). Now, we add these numbers using the relation (\ref{1.3}\ (2)) 
 applied for $mod_N(G)({\bf Str})$\
(namely, for any $w\in\calw$, the relation 
$e_wm_w+\sum_{v\in\calv_w}m_v=0$ provides $e_w$).

\bekezdes{Theorem.}{7.2} {\em 
If $\phi:({\caly},D)\to(X,x)$ is an embedded resolution
of 
$(f^{-1}(0),x)\subset (X,x)$, and 
$p:G(X,f)\to \G(X,f)$ the  universal covering graph 
associated with $\phi$. Then
the graphs $\G(X_{f,N},z)$ and $\G(X_{f,N})$ can be determined from 
$p:G(X,f)\to \G(X,f)$ and from the integer $N$.
Namely:

a)\ The decorated graph $mod_N(G)({\bf Str})$ constructed in (\ref{7.1})
is (a possible) embedded resolution graph $\G(X_{f,N},z)$.

b)\ If we delete all the arrows and multiplicities of $mod_N(G)({\bf Str})$, 
then we obtain a resolution graph $\G(X_{f,N})$ of $X_{f,N}$.

}

\bekezdes{Example.}{7.3} 
Set $(X,x)=(\{x^2+y^7-z^{14}=0\},0)\subset
(\bfc^3,0)$. Take $f_1(x,y,z)=z^2+y^2$ and $f_2(x,y,z)=z^2-y+y^2$ 
(cf. \ref{3.17} and \ref{5.10}).
Then $\G(X,f_1)=\G(X,f_2)$, but in general, the graphs $\G(X_{f_i,N},z)$ 
($i=1,2$) are 
not the same. For example, these graphs for $N=2$ are:

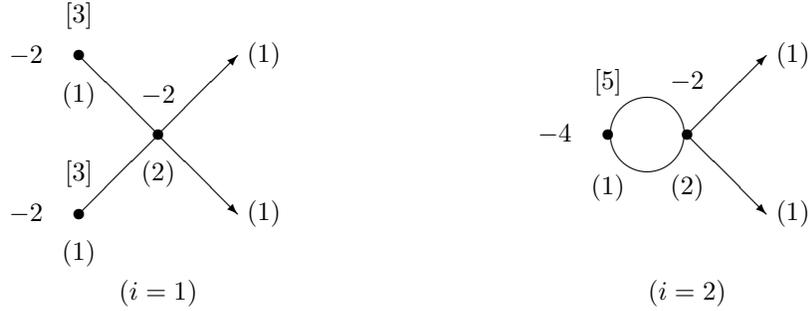
\begin{picture}(350,140)(-20,50)

\put(10,165){\makebox(0,0){$[3]$}}
\put(10,135){\makebox(0,0){$(1)$}}
\put(10,105){\makebox(0,0){$[3]$}}
\put(10,75){\makebox(0,0){$(1)$}}

\put(-10,150){\makebox(0,0){$-2$}}
\put(-10,90){\makebox(0,0){$-2$}}

\put(40,135){\makebox(0,0){$-2$}}
\put(40,105){\makebox(0,0){$(2)$}}
\put(80,150){\makebox(0,0){$(1)$}}
\put(80,90){\makebox(0,0){$(1)$}}

\put(40,60){\makebox(0,0){$(i=1)$}}
\put(240,60){\makebox(0,0){$(i=2)$}}

\put(10,90){\circle*{4}}
\put(10,150){\circle*{4}}
\put(40,120){\circle*{4}}
\put(10,90){\line(1,1){30}}
\put(10,150){\line(1,-1){30}}
\put(40,120){\vector(1,1){30}}
\put(40,120){\vector(1,-1){30}}

\put(210,140){\makebox(0,0){$[5]$}}
\put(210,100){\makebox(0,0){$(1)$}}
\put(190,120){\makebox(0,0){$-4$}}

\put(240,140){\makebox(0,0){$-2$}}
\put(240,100){\makebox(0,0){$(2)$}}
\put(280,150){\makebox(0,0){$(1)$}}
\put(280,90){\makebox(0,0){$(1)$}}

\put(210,120){\circle*{4}}
\put(240,120){\circle*{4}}
\put(225,120){\circle{29}}
\put(240,120){\vector(1,1){30}}
\put(240,120){\vector(1,-1){30}}

\end{picture}

\noindent 
For $N=4$, the graphs $\G(X_{f_i,N},z)$ are:

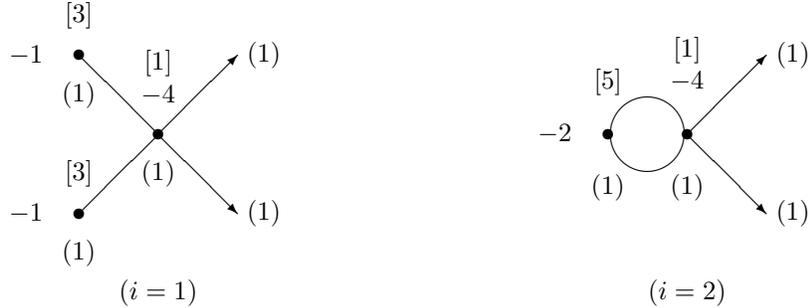
\begin{picture}(350,140)(-20,50)

\put(10,165){\makebox(0,0){$[3]$}}
\put(10,135){\makebox(0,0){$(1)$}}
\put(10,105){\makebox(0,0){$[3]$}}
\put(10,75){\makebox(0,0){$(1)$}}

\put(-10,150){\makebox(0,0){$-1$}}
\put(-10,90){\makebox(0,0){$-1$}}

\put(40,147){\makebox(0,0){$[1]$}}
\put(40,135){\makebox(0,0){$-4$}}
\put(40,105){\makebox(0,0){$(1)$}}
\put(80,150){\makebox(0,0){$(1)$}}
\put(80,90){\makebox(0,0){$(1)$}}

\put(40,60){\makebox(0,0){$(i=1)$}}
\put(240,60){\makebox(0,0){$(i=2)$}}

\put(10,90){\circle*{4}}
\put(10,150){\circle*{4}}
\put(40,120){\circle*{4}}
\put(10,90){\line(1,1){30}}
\put(10,150){\line(1,-1){30}}
\put(40,120){\vector(1,1){30}}
\put(40,120){\vector(1,-1){30}}

\put(210,140){\makebox(0,0){$[5]$}}
\put(210,100){\makebox(0,0){$(1)$}}
\put(190,120){\makebox(0,0){$-2$}}

\put(240,152){\makebox(0,0){$[1]$}}
\put(240,140){\makebox(0,0){$-4$}}
\put(240,100){\makebox(0,0){$(1)$}}
\put(280,150){\makebox(0,0){$(1)$}}
\put(280,90){\makebox(0,0){$(1)$}}

\put(210,120){\circle*{4}}
\put(240,120){\circle*{4}}
\put(225,120){\circle{29}}
\put(240,120){\vector(1,1){30}}
\put(240,120){\vector(1,-1){30}}

\end{picture}

\noindent 
For any odd $N$ the orbit graphs $mod_N(G)$, for $i=1,2$, are the same. 
Hence  in this case
$\G(X_{f_1,N},z)=\G(X_{f_2,N},z)$.   For $N=3$ this graph is:

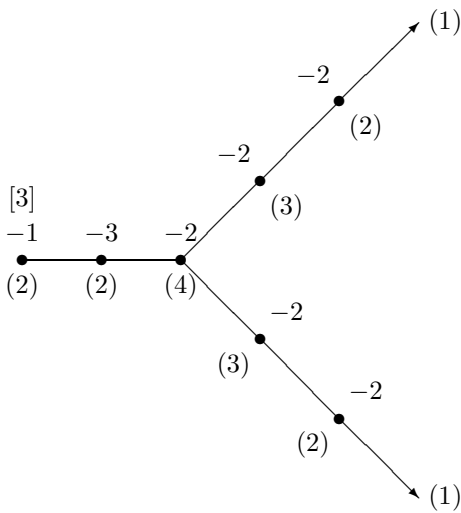
\begin{picture}(300,200)(0,-80)
\put(80,20){\circle*{4}}
\put(110,20){\circle*{4}}
\put(140,20){\circle*{4}}
\put(140,20){\vector(1,1){90}}
\put(140,20){\vector(1,-1){90}}

\put(170,50){\circle*{4}}
\put(200,80){\circle*{4}}
\put(170,-10){\circle*{4}}
\put(200,-40){\circle*{4}}

\put(160,60){\makebox(0,0){$-2$}}
\put(190,90){\makebox(0,0){$-2$}}
\put(180,0){\makebox(0,0){$-2$}}
\put(210,-30){\makebox(0,0){$-2$}}

\put(180,40){\makebox(0,0){$(3)$}}
\put(210,70){\makebox(0,0){$(2)$}}
\put(160,-20){\makebox(0,0){$(3)$}}
\put(190,-50){\makebox(0,0){$(2)$}}

\put(110,30){\makebox(0,0){$-3$}}
\put(80,30){\makebox(0,0){$-1$}}
\put(140,30){\makebox(0,0){$-2$}}
\put(80,43){\makebox(0,0){$[3]$}}
\put(80,20){\line(1,0){60}}

\put(110,10){\makebox(0,0){$(2)$}}
\put(80,10){\makebox(0,0){$(2)$}}
\put(140,10){\makebox(0,0){$(4)$}}
\put(240,110){\makebox(0,0){$(1)$}}
\put(240,-70){\makebox(0,0){$(1)$}}
\end{picture}

\bekezdes{Example.}{7.4}\ 
Set $(X,x)=(\{z^2+(x^2-y^3)(x^3-y^2)=0\},0)$
and $f_1=x^2+y^4$ and $f_2=x^2-y^3+y^4$ (cf. \ref{3.18} and \ref{5.12}).
The graphs $\G(X_{f_i,N},z)$  for $N=2$ are:

\hspace{-1.3cm}
\begin{picture}(400,140)(20,50)
\put(110,70){\makebox(0,0){$(i=1)$}}

\put(80,120){\circle*{4}}
\put(110,110){\circle*{4}}
\put(110,130){\circle*{4}}
\put(140,160){\circle*{4}}
\put(140,140){\circle*{4}}
\put(170,100){\circle*{4}}
\put(170,80){\circle*{4}}
\put(80,120){\line(3,1){30}}
\put(80,120){\line(3,-1){30}}
\put(110,110){\line(1,1){30}}
\put(110,130){\line(1,1){30}}
\put(110,110){\line(2,-1){60}}
\put(110,130){\line(2,-1){60}}
\put(140,160){\line(1,-2){30}}
\put(140,140){\line(1,-2){30}}
\put(80,120){\vector(-1,-1){30}}
\put(80,120){\vector(-1,1){30}}

\put(80,110){\makebox(0,0){$a$}}
\put(110,140){\makebox(0,0){$b$}}
\put(110,100){\makebox(0,0){$b$}}

\put(280,110){\makebox(0,0){$a$}}
\put(310,140){\makebox(0,0){$b$}}
\put(310,100){\makebox(0,0){$b$}}

\put(310,70){\makebox(0,0){$(i=2)$}}

\put(280,120){\circle*{4}}
\put(310,110){\circle*{4}}
\put(310,130){\circle*{4}}
\put(340,160){\circle*{4}}
\put(340,140){\circle*{4}}
\put(370,100){\circle*{4}}
\put(370,80){\circle*{4}}
\put(280,120){\line(3,1){30}}
\put(280,120){\line(3,-1){30}}
\put(310,110){\line(3,5){30}}

\put(310,130){\line(3,1){30}}
\put(310,110){\line(2,-1){60}}
\put(310,130){\line(2,-1){60}}
\put(340,160){\line(1,-2){30}}
\put(340,140){\line(1,-2){30}}
\put(280,120){\vector(-1,-1){30}}
\put(280,120){\vector(-1,1){30}}

\end{picture}

The decorations are the following. The arrows have multiplicity $(1)$,
$m_a=4,\ m_b=3$, and the multiplicities of the unmarked nodes 
are $m=1$; $e_a=e_b=-2$, and the self intersection of the unmarked
nodes are $e=-4$.  All the genera are zero.

\bekezdes{Remark.}{7.5}\ The construction (\ref{7.1}) gives 
non--minimal resolution graphs, in general. They can be simplified by blowing--down
the $(-1)$--rational exceptional divisors $E_w$ with $\delta_w\leq 2$.

\bekezdes{Remark.}{7.6}\ Using the formula of (Step 2.), it is easy to prove
 that $\tilde{g}_{\tilde{w}}\geq g_w$. Moreover,  $c_{mod_N(G)}\geq c_{\G}$
(cf. also with (\ref{cgcg}) and
(\ref{9.8})). Therefore:
$$rank\ H_1(L_{X_{f,N}})\geq rank\ H_1(L_X).$$

\newpage

\bekezdes{-- Modification of (\ref{7.1}) for the case when $L_X$ is a rational
homology sphere.}{8.1}
If $L_X$ is a rational homology sphere, then the above algorithm can be
simplified. In this case the universal covering $\gxf\longrightarrow\Gxf$
can completely be reconstructed from $\Gxf$ (see \ref{5.5}--\ref{L}).
This means that the graph $\Gamma (\xfn, z)$ is completely determined by
the embedded resolution graph $\Gxf$ of $f$ and the integer $N$, and the
reader can easily reconstruct this new algorithm. 

In particular, if $(X,x)=(\bfc ^2,0)$, and 
$f: (\bfc ^2,0)\longrightarrow (\bfc,0)$ is an arbitrary 
isolated plane curve singularity, then 
$(\xfn ,0)=(\fegyz ,0)\subset (\bfc^3,0)$, and $z:(\xfn ,0)\longrightarrow
(\bfc,0)$ is induced by the projection $(x,y,z)\longmapsto z$.
Therefore, our algorithm also provides the resolution graph of 
$(Y,0)=(\{ g=0\} ,0)$, where $g=f(x,y)+z^N$.
For the algorithm in this case  see  \cite{NRev}. The idea of the
construction in the case of $f+z^N$ can already 
be found in the book of Laufer \cite{La}.
For other particular cases, see \cite{Artal,O,OW}. 

\vspace{1mm}

The general algorithm  can be compared with some results of Eriko Hironaka, who
considers the global case of cyclic coverings. 

\bekezdes{Proof of the Theorem (\ref{7.2}).}{proof} \ 
Consider the ``Jungian diagram'' (cf. \ref{2.1}):

\begin{picture}(400,100)(0,0)
\put(50,75){\makebox(0,0){$X^{res}$}}
\put(100,75){\makebox(0,0){$\tilde{X}'$}}
\put(160,75){\makebox(0,0){$X'$}}
\put(230,75){\makebox(0,0)[l]{$(X_{f,N},x)$}}
\put(160,25){\makebox(0,0){$({\caly},D)$}}
\put(230,25){\makebox(0,0)[l]{$(X,x)$}}
\put(65,75){\vector(1,0){20}}
\put(115,75){\vector(1,0){20}}
\put(190,75){\vector(1,0){30}}
\put(190,25){\vector(1,0){30}}
\put(115,60){\vector(1,-1){20}}
\put(160,60){\vector(0,-1){20}}
\put(250,60){\vector(0,-1){20}}
\put(60,60){\vector(2,-1){50}}
\put(60,50){\makebox(0,0){$\rho$}}
\put(125,85){\makebox(0,0){$n$}}
\put(170,50){\makebox(0,0){$\pi'$}}
\put(115,50){\makebox(0,0){$\tilde{\pi}'$}}
\put(260,50){\makebox(0,0)[l]{$\pi=pr$}}
\put(200,85){\makebox(0,0){$\phi'$}}
\put(200,35){\makebox(0,0){$\phi$}}
\end{picture}

\noindent where:

a)\ the proper map
$\pi:(X_{f,N},x)\to (X,x)$ is induced by the projection $(x',z)\mapsto
x'$, and it is an $N$--covering with branch locus $f^{-1}(0)$ (cf. 
\ref{6.1}).

b)\ $\phi$ is a fixed  embedded resolution of $(f^{-1}(0),x)\subset
(X,x)$, and  $D$ is the divisor $\phi^{-1}(f^{-1}(0)) $. The dual graph
of $D$ is exactly $\G(X,f)$.

c)\ $\pi':X'\to {\caly}$ is the pullback of $\pi$ via $\phi$, and
$\tilde{X}'$ is the normalization of $X'$. Then $\tilde{X}'$ has only 
Hirzebruch--Jung singularities (cf. \ref{1.9} and \ref{2.1}).

d)\ $X^{res}\to \tilde{X}'$ is the resolution of the Hirzebruch--Jung
singularities of $\tilde{X}'$.\\

If $P$ is a generic point of the irreducible exceptional divisor 
$E_w\subset E=\phi^{-1}(x)$ ($w\in\calw(\G(X,f))$), then we fix local
coordinates $(u,v)$ in a neighbourhood $U$ of $P$ such that $\{u=0\}=
E_w\cap U$, and $f\circ \phi|_{U}=u^{m_w}$. Let $P'$ be the unique point in 
$X'$ above $P$ and consider  its neighbourhood $U':=(\pi')^{-1}(U)$.
Then $U'=\{(u,v,z)\in (U\times \bfc,P\times 0): u^{m_w}=z^N\}$.
This shows that $n^{-1}(P')$ contains exactly $(N,m_w)$ points, which 
correspond to the irreducible components of $U'$. Indeed, $U'=\cup_i
U'_i$, where $U'_i=\{(u,v,z): u^{m_w/(N,m_w)}=\epsilon_i\cdot z^{N/(N,m_w)}\}$,
 where $\epsilon_i$, for $i=0,1,\ldots, (N,m_w)-1$,  are the $(N,m_w)$-roots 
of unity. Moreover, the normalization of any $U'_i$ is smooth. For example,
for $i=0$, a possible normalization  map of $U'_0$ is $(\bfc^2,0)\to
U'_0$ given by $(t,v)\mapsto(u(t),v,z(t))$, where 
$$(*)\hspace{2cm}
u(t)=t^{N/(N,m_w)}, \ \ z(t)=t^{m_w/(N,m_w)}.\hspace{2cm}$$
Therefore, $(\tilde{\pi}')^{-1}(E_w)\to E_w$ is a regular covering of degree
$(N,m_w)$ above $E_w^0:=E_w\setminus\overline{D\setminus E_w}$. The covering 
$U'\to U$ is a Galois covering with Galois group action of
$\bfz_{N}=\{\xi:\xi^N=1\}$ given by $(u,v,z)\mapsto (u,v,\xi z)$. Notice that 
$\xi^{m_w}$ preserves the components $U'_i$, therefore the covering
$(\tilde{\pi}')^{-1}(E_w^0)\to E_w^0$ is a Galois covering with Galois group
$\bfz_N/m_w\bfz_N=\bfz_{(N,m_w)}$. This regular covering extends to a branched
covering $(\tilde{\pi}')^{-1}(E_w)\to E_w$, with branch locus $E_w\setminus E_w^0$. The curve 
$(\tilde{\pi}')^{-1}(E_w)$ contains exactly the irreducible exceptional
 divisors of $\tilde{X}'$ lying above $E_w$. Then the number $i_w$ of these 
irreducible exceptional 
components is exactly the number of irreducible components
of $(\tilde{\pi}')^{-1}(E_w^0)$, and this number coincides with the number
of connected components of $(\tilde{\pi}')^{-1}(E_w^0)$. The representation
associated with the covering $(\tilde{\pi}')^{-1}(E_w^0)\to E_w^0$ is
denoted by:
$$\rho_w:\pi_1(E_w^0)\to\bfz_{(N,m_w)},$$
and the number $i_w$ of connected components of $(\tilde{\pi}')^{-1}(E_w^0)$
is $|\coker(\rho_w)|$. 

Now, let $T_w$ be a tubular neighborhood of $E_w$ in ${\caly}$ (cf. \ref{5.1}).
Then $(\tilde{\pi}')^{-1}(T_w\setminus D)\to T_w\setminus D$ is a regular 
$N$--covering. The corresponding representation is denoted by:
$$r_w: \pi_1(T_w\setminus D) \to\bfz_N.$$ 
This representation is induced by $arg_*$, as it is explained in
(\ref{6.1}--\ref{6.2}). Hence, the following composed map:
$$\pi_1(T_w\setminus D)\stackrel{j}{\longrightarrow}\pi_1(L_X\setminus L_f)
\stackrel{arg_*}{\longrightarrow}\bfz\stackrel{pr_N}{\longrightarrow}\bfz_N$$
is exactly the map $r_w$ (above, $j$ is induced by the natural inclusion).

Now, notice that $T_w\setminus D=E_w^0\times S^1$, where $S^1$ can be represented 
by an oriented circle in a transversal slice of $E_w^0$ in $T_w$. In particular,
$\pi_1(T_w\setminus D)=\pi_1(E_w^0)\times \bfz$, and $arg_*((0,1_{\bfz}))=m_w\in\bfz$.
 This shows that $r_w((0,1_{\bfz}))$ is the class of $m_w$ in $\bfz_N$. The above discussion
provides the following commutative diagram:

$$\begin{array}{rccc}
\pi_1(T_w\setminus D)=&\pi_1(E_w^0)\times \bfz &\stackrel{r_w}{\longrightarrow}&\bfz_N\\
 &\uparrow & & \downarrow\vcenter{%
\rlap{$\scriptstyle{pr_{(N,m_w)}}$}}\\
& \pi_1(E_w^0)&\stackrel{\rho_w}{\longrightarrow}&\bfz_{(N,m_w)}\end{array}$$

Since the class of $m_w$ in $\bfz_{(N,m_w)}$ is zero, one has:
$|\coker(pr_{(N,m_w)}\circ r_w)|=|\coker(\rho_w)|$. 

Therefore, $i_w=|\coker(\rho_w)|$ is the cardinality of the cokernel of the 
following composed map:
$$p_w:\pi_1(T_w\setminus D)\stackrel{arg_*\circ j}{\longrightarrow}\bfz
\stackrel{pr_N}{\longrightarrow}\bfz_N\stackrel{pr_{(N,m_w)}}{\longrightarrow}
\bfz_{(N,m_w)}.$$
Now, consider the covering $p:G(X,f)\to \G(X,f)$, and, as usual, let $n_w$ be
 the number of connected components of the Milnor fiber $F$ in $T_w$ (i.e.
the number of vertices of $G$ above $w\in\calw(\G)$). But this is exactly
$|\coker(arg_*\circ j)|$, in other words: $im(arg_*\circ j)=n_w\bfz$. 
Since $arg_*((0,1))=m_w\in\bfz$, one has: $n_w|m_w$. Hence $|\coker(p_w)|=
(N,n_w)$. But this is exactly the number of vertices above $w$ in $mod_N(G)$.\\

The above argument can be repeated in the case of the singular points of the 
divisor $D$. Let $P$ be  an intersection point of two irreducible components 
of $D$, and fix local coordinates $(u,v)$ at $P$ such that the local equation of $D$
at $P$ is $\{uv=0\}$, and $\phi\circ f$ in a neighborhood of $P$ is $u^{m_w}v^{m_v}$.
Set $P'=(\pi')^{-1}(P)$. Then the local equation of $X'$ at $P'$ is $\{
u^{m_w}v^{m_v}=z^N\}$, which has $(N,m_w,m_v)$ irreducible components, hence
$(\tilde{\pi}')^{-1}(P)$ contains exactly $(N,m_w,m_v)$ points. Since $n_e=
(m_w,m_v)$ is the number of edges of $G$ above the edge $(v,w)\in\cale(\G)$,
\ $(N,n_e)=(N,m_w,m_v)$ is the number of edges in $mod_N(G)$ above $(v,w)$.
Hence it coincides with $\#(\tilde{\pi}')^{-1}(P)$.  
Moreover, the following diagram provides the adjacency relations in $G$:

$$\begin{array}{cccccl}
\pi_1(T_w\setminus D)&\stackrel{arg_*\circ j}{\longrightarrow}
&\bfz&\to &\bfz_{n_v}&\to 0\\
\uparrow& & \uparrow \vcenter{\rlap{$=$}}& & 
\uparrow \vcenter{\rlap{$\scriptstyle{a}$}}& \\
\pi_1(T_w\cap T_v \setminus D)&\stackrel{arg_*\circ j}{\longrightarrow}
&\bfz&\to &\bfz_{n_e}&\to 0
\end{array}$$

The group $\bfz_{n_w}$  (resp. $\bfz_{n_e}$) is the index set of the vertices 
(resp. edges) of $G$ above $w$ (resp. above $e=(w,v)$). The edge in $G$ indexed by
$\hat{l}\in\bfz_{n_e}$ has as one of its endpoints the vertex indexed by
$a(\hat{l})$. 

Now, in the commutative diagram:

$$\begin{array}{ccc}
\bfz_{n_w}&\to& \bfz_{(N,n_w)}\\
\uparrow\vcenter{\rlap{$\scriptstyle{a}$}} & & \uparrow\vcenter{
\rlap{$\scriptstyle{a}$}}\\
\bfz_{n_e}&\to& \bfz_{(N,n_e)}
\end{array}$$

\noindent the left column codifies the adjacency relations in $G$, while the
right column the adjacency relations in $mod_N(G)$ (in  a similar way). 
Therefore, by 
the above discussion, it is clear that $mod_N(G)$ is the dual graph of the 
exceptional divisors of $\tilde{X}'$.\\

The space $\tilde{X}'$ has only Hirzebruch--Jung singularities. The resolution
of these singularities does  not change the multiplicities of $z$ along the 
irreducible exceptional divisors of $\tilde{X}'$, nor the topology of 
them. So, already at the level of $\tilde{X}'$ we can compute these invariants.
First notice that the irreducible components above $E_w$ are cyclically 
permuted by the Galois action, so their multiplicities and genera are the same.
The formula $(*)$ (of this proof)
shows that at a generic point of $(\tilde{\pi}')^{-1}(E_w)$,
the multiplicity of $z$ is $m_w/(N,m_w)$, proving (Step 2, a). 

Above the branch point $E_v\cap E_w$ of the 
covering $ (\tilde{\pi}')^{-1}(E_w)\to E_w$, there are $(N,m_w,m_v)$ points,
therefore, by an Euler--characteristic argument:

$$\chi((\tilde{\pi}')^{-1}(E_w))=(2-2g_w-\delta_w)(N,m_w)+\sum_{v\in\calv_w}
(N,m_w,m_v).$$
But this is $(N,n_w)\cdot (2-2\tilde{g})$, proving (Step 2,b).\\

Now, we have to resolve the Hirzebruch--Jung singularities
of  $\tilde{X}'$. The singular points  $\{P_{\tilde{e}}\}$ 
of  $\tilde{X}'$ are codified by the 
edges $\tilde{e}$ of $mod_N(G)$, and their local equations are
$\{u^{m_w/t_e}v^{m_v/t_e}=z^{N/t_e}\}$, where $t_e=(N,m_w,m_v)$, \ $e=(v,w)$.
Therefore, any edge $\tilde{e}$ of $mod_N(G)$ must be replaced by the 
string of the corresponding Hirzebruch--Jung singularity. Hence the last part 
follows from subsection (\ref{2.3}).\  \ \  \hfill $\diamondsuit$

\subsection*{Cyclic coverings and monodromy.}

\vspace{2mm}

\bekezdes{}{9.1} The main goal of the present section is to compare the
covering $p:G(X,f)\to \G(X,f)$ with the algebraic monodromy operator
associated with the Milnor fibration $arg:L_X\setminus L_f\to S^1$.

Start again with a normal surface singularity $(X,x)$ and a germ
$f:(X,x)\longrightarrow(\bfc,0)$ of an analytic function. 
The characteristic map of the Milnor fibration 
$F\longrightarrow L_X\setminus L_f\stackrel{arg}\longrightarrow S^1$ is called the 
geometric monodromy $\geo :F\longrightarrow F$ (where $F$ denotes the Milnor 
fiber). 

At homology level $\geo$ induces the algebraic monodromies
$h_q:H^q(F,\bfz )$ $\longrightarrow H^q(F,\bfz )$  \ $ (q=0,1)$.
If $F$ has $k$ connected components, then $H^0(F,\bfz )\approx\bfz^k$ and 
$h_0(x_1,...,x_k)=(x_2,...,x_k,x_1)$ (cf. \ref{3.12}).
In particular, $h_0$ is finite.

The monodromy $h_1$ is more complicated -- in general it is not finite.
But, by the Monodromy Theorem (see, e.g. 
\cite{BrieskornMon,Landman,LeMon,A'CMon,A'CZeta}, 
all the eigenvalues of $h_1$ are 
roots of unity, and its Jordan decomposition contains blocks only 
of size one and two. 

In this section, we connect the Jordan decomposition of $h_1$
with the number of 
independent cycles in the cyclic coverings with branch locus $\{f=0\}$.
Finally, we show that the number of Jordan blocks of $h_1$
is completely determined by the universal cyclic covering graph
$\gxf\longrightarrow\Gxf$. Some results of this subsection can be
compared with some results of Ph. Du Bois and F. Michel \cite{MDB}. 

\bekezdes{-- The characteristic polynomial.}{9.2}
 For $q=0,1$ we define $\Delta_q(t)=\det (tI- h_q)$.
Then, for $q=0$, $\Delta_0(t)=t^k-1$. The characteristic polynomial
$\Delta_1(t)$ is determined by $\Gxf$, via A'Campo's formula 
\cite{A'CMon,A'CZeta}:

\bekezdes{}{9.3}\ \

\vspace{-1cm}

\
$$
\Delta_0(t)/\Delta_1(t) =\prod\limits_{w\in\calw (\G )} 
(t^{m_w}-1)^{2-2g_w-\delta_w}.$$

\noindent 
It is convenient to use the notation $H=H^1(F,\bfc )$; and 
 $H_\lambda$ for the generalized eigenspace corresponding to the
eigenvalue $\lambda$ of $h_1$ 
(i.e. $H_\lambda=\{ x\in H \vert\ (h_1-\lambda I)^nx=0$
for some sufficiently\ large \ $n\}$.)  
The dimension $\dim_\bfc H_\lambda$ is 
exactly the order of $t-\lambda$ in $\Delta_1(t)$.

\bekezdes{-- The Jordan blocks of $h_1$.}{9.4}

Let  $\szam{\lambda}{l}$ be the number of Jordan 
blocks of size $l$ of the 
restriction $h_1\vert_{H_\lambda}$.
Obviously,  $\szam{\lambda}{1}+2\szam{\lambda}{2}=\dim H_\lambda$
and $\dim \ker (h_1-\lambda I)=
\szam{\lambda}{1}+\szam{\lambda}{2}.$\\

\noindent First consider $\lambda =1$.

\bekezdes{Proposition.}{9.5} (cf. \cite{NSM}) {\em 
Let $c_{\G (X)}$ be the number of independent cycles in  $\G (X)$, and 
$g=\sum_{w\in\calw(\G)}g_w$. Then:

a)\   $\dim H_1=2g+2c_{\G (X)}+\#\cala (\G )-1; $

b)\   $\dim\ker (h_1-1)=2g+c_{\G (X)}+\#\cala (\G )-1. $

Therefore, $\szam{1}{2}=c_{\G(X)}$. 
In particular, $\szam{1}{2}$ is independent of the germ $f$, and depends
only on the topology of the link $L_X$.}

\vspace{2mm}

\noindent {\em Proof.}\  The order of $t-1$ in
$\Delta_0/\Delta_1$ is $1-\dim H_1$. Via A'Campo's formula this is
$\sum_w (2-2g_w-\delta_w)$.
Now, use $\sum_w\delta_w =\#\cala +2\#\cale$, and the 
``Euler-characteristic'' identity $1-c_{\G (X)}=\#\calw -\#\cale $, and $(a)$
follows. By the Wang exact sequence associated with the Milnor fibration
\[
H^0(F)\stackrel{h_0-I}\longrightarrow H^0(F)\longrightarrow H^1(L_X\setminus L_f)
\longrightarrow H^1(F)\stackrel{h_1-I}\longrightarrow H^1(F)
\]
one has that $\coker(h_0-I)\approx\bfc$, hence $\dim\ker (h_1-1)=\dim
H^1(L_X\setminus L_f)-1$. But using (\ref{3.8}) and (\ref{3.9})
$\dim H^1(L_X\setminus L_f)=\dim H_1(E)+\#\cala =2g+c_\G +\#\cala$.
\hfill $\diamondsuit$

\vspace{2mm}

As a corollary of (\ref{9.5}), we obtain that if $L_X$ has no cycles (e.g. 
$(X,x)$ is smooth as in the case of plane curve singularities), then $f$ has no
Jordan block of size 2 with eigenvalue $\lambda =1$.

\smallskip
Now, assume that $f:(X,x)\longrightarrow(\bfc,0)$ is the smoothing of 
 an {\em isolated singularity}
at $x$. Then $arg_*(f)$ is onto, hence for any $N$, $\xfn$ is connected.
Moreover, 
there is  natural identification of the Milnor fibers  of $f$ and $z:
(X_{f,N},x)\to (\bfc,0)$ such that 
the monodromy of $z$  is the 
$N^{th}$--power of the monodromy of $f$. Therefore,  by (\ref{9.5}) 
$c_{\G (\xfn)} =\#_{1}^{2}(z)=\sum\limits_{\lambda^N=1}\szam{\lambda}{2}$.
Thus we have the following:

\bekezdes{Corollary.}{9.6}\
{\em  If $f$ defines an isolated singularity, then:}
\[
\sum\limits_{\lambda^N=1}\szam{\lambda}{2} = c_{\G (\xfn)}.
\]

Since any eigenvalue of $h_1$ is a root of unity,
$h_1$ is of finite order if and only if it has no Jordan blocks of size 2.
Therefore, the above corollary 
generalizes a theorem of A. Durfee, which says that for 
plane curve singularities $f:(\bfc^2,0)\longrightarrow(\bfc,0)$, 
$h$ has finite order
  if and only if the graphs of the cyclic coverings have no cycles 
\cite{DuF}.

\bekezdes{Remark.}{9.7}\ 
In some particular cases $c_{\G (\xfn)}$ can be computed from the 
multiplicity system of $\Gxf$ alone. For example, if $L_X$ is a rational 
homology sphere, then by (\ref{5.5}):
\[
\begin{array}{rcl}
    1-c_{\G (\xfn)} &=&\#\calw (\G (\xfn ))-\#\cale (\G (\xfn ))\\ 
                    & &       \\
                    &=&\sum\limits_{w\in\calw(\Gxf )} (M_w,N) -
                       \sum\limits_{e\in\cale_n(\Gxf )} (m_e,N) 
\end{array}
\]

\noindent where $M_w=g.c.d. \{ \{ m_v\} _{v\in\calv_w}, m_w \}$ and 
$m_e=g.c.d.(m_{v_1},
m_{v_2})$ with $e=(v_1,v_2)$, and $\cale_n$ denotes the set of edges
connecting two non--arrowheads. Then by (\ref{9.6}) and by an easy computation
using $\#\calw-\#\cale_n=1$ one has:
\[
\sum\limits_{\lambda^N=1}\szam{\lambda}{2}=\sum\limits_{e\in\cale_n
 (\Gxf )} ((m_e,N)-1)-\sum\limits_{w\in\calw (\Gxf )}((M_w,N)-1)
\]

This is nothing else than W. Neumann's formula for $\#_{\lambda}^2(f)$
\cite{EN,Ne2}.\\

Now,  we will show that, in order to determine the numbers
$\#^2_{\lambda}(f)$,  we don't have to 
consider all the cyclic coverings, but only 
$\gxf$, the universal cyclic covering graph of $\Gxf$.
First notice that the identity 
(\ref{9.6}) is valid even if $\finv$ has non-isolated 
singularities, but $\gxf$ is connected. Then $\xfn$ is connected for any $N$.

\bekezdes{Corollary.}{9.8}\
{\em  Assume that the universal covering graph $\gxf$ is 
connected  (e.g. \  $f$ defines an isolated singularity).
Then:

(a)$$ \sum_{\lambda^N=1}\#_{\lambda}^2(f)=c_{mod_N(G(X,f))}.$$

(b)$$  \sum_{\lambda}\#_{\lambda}^2(f)=c_{G(X,f)};$$
i.e. the number of independent cycles of $\gxf$ is exactly the total number of
Jordan blocks of size 2 of the algebraic monodromy $h_1$ of $f$.

(c)\ Let $\xi_n$ be a primitive $n^{th}$--root of unity. Then:
$$\phi(n)\cdot \#^2_{\xi_n}=\sum_{k|n}\mu(n/k)\cdot c_{mod_k(G(X,f))}$$
where $\phi$ is the Euler function and $\mu$ is the M\"obius function,
namely: $\phi(n)=\#\{1\leq k\leq n:\ (k,n)=1\}$, and $\mu(k)=1$ for $k=1$,
$\mu(k)=(-1)^t$ if $k$ is a product of $t$ different primes, and $\mu(k)=0$ if 
$p^2|k$ for some prime $p$.
}

\vspace{2mm}

\noindent {\em Proof.}\ (a) follows from a similar argument as (\ref{9.6}) 
and the main theorem (\ref{7.2}).
For (b), take 
an $N$ which is multiple of all the integers 
$\{ n_v\}_{v\in\calv}$ and $\{ n_e\}_{e\in\cale}$ (the integers which
characterize the type of 
 the universal covering $\gxf\longrightarrow\Gxf$), and also $\lambda^N=1$
for any eigenvalue $\lambda$ of $h_1$. For (c), first notice that 
$\#^2_{\xi_n}$ does not depend on the choice of $\xi_n$
(because $h_1$ is defined over $\bfz$). Then (c)
 follows from the M\"obius Inversion
Formula (see, e.g. \cite{Hua}, page 107).\hfill  $\diamondsuit$

\bekezdes{Example.}{9.9} Set $(X,x)$ as in (\ref{5.10}) (cf. also
(\ref{3.17})). 
Then $h_1(f_1)$ is finite, but $h_1(f_2)$ has a Jordan block 
of size 
2, with eigenvalue $\lambda =-1$. In this case $\G (X,f_1)=\G (X,f_2)$, but
$arg_*(f_1)\neq arg_*(f_2)$. This shows a subtle connection between the
Jordan block structure of $h_1(f)$ and the representation $arg_*(f)$.

\bekezdes{Example.}{9.10}
Set $(X,x)$ as in (\ref{5.12}) (or \ref{3.18}).
Then for both $i=1,2$, $h_1(f_i)$ has a Jordan block of size 2 with eigenvalue
$\lambda =1$, because $\Gxf$ has a cycle. Since $\gxf$ has two independent
cycles, there is one more Jordan block of size 2 and this one has eigenvalue
$\lambda =-1$.

\bekezdes{}{9.11}
Actually, an even stronger connection can be established between the monodromy
$h_1$ and the graph $G(X,f)$. Let $|G|$ be the topological realization of 
$G(X,f)$. The $\bfz$--action of $G(X,f)$ induces a ``geometric action''
$h_{G,geom}$ on $|G|$ (by the identification $1_{\bfz}=
h_{G,geom}$). At homological level, this induces a finite morphism 
$h^*_{|G|}:H^1(|G|)\to H^1(|G|)$. 

In the sequel, $H^*(Y)$ denotes the cohomology with complex coefficients.
We invite the reader to review the notations and the results of 
subsection (\ref{5.1}).
Recall that $\calv(\G)$ is the index set of the irreducible components of
$D=\phi^{-1}(f^{-1}(0))$. Let $T_v$ be a small tubular neighborhood
of the irreducible component corresponding to $v\in\calv(\G)$. For any edge
$e=(v_1,v_2)$ of $\G$ set $T_e=T_{v_1}\cap T_{v_2}$. If $F=f^{-1}(\delta)$
is the Milnor fiber, then for $\delta$ sufficiently small $F\subset 
\cup_vT_v$. Put $F_v=F\cap T_v$ and $F_e=F\cap T_e$ for any $v\in\calv(\G)$
and $e\in\cale(\G)$. 
Then  by Mayer--Vietoris argument one has the following exact sequence:
$$0\to H^0(F)\to\oplus_vH^0(F_v)\to\oplus_eH^0(F_e)\to
H^1(F)\to\oplus_vH^1(F_v)\to\oplus_eH^1(F_e)\to 0.$$

\noindent Now, by the very definition 
of the graph $G=G(X,f)$, one has the exact sequence:
$$0\to H^0(|G|)\to\oplus_vH^0(F_v)\to\oplus_eH^0(F_e)\to H^1(|G|)\to 0.$$

Now, we recall that there is a continuous map $\alpha:F\to |G|$ 
(cf. \ref{5.1}) which is compatible with the above decomposition, in particular
the above exact sequences (complexes)  are connected by maps 
(i.e. by a morphism of complexes) induced by  $\alpha$. 

This  provides an exact sequence:
$$0\to H^1(|G|)\stackrel{\alpha^*}{\to}H^1(F)\stackrel{\delta^0}{\to}
\oplus_vH^1(F_v)\stackrel{\delta^1}{\to}\oplus_eH^1(F_e)\to 0.$$

The monodromy acts on this exact sequence, the operators on the corresponding
groups are: $h^*_{|G|},\ h_1,\ \oplus_wh_w$ and $\oplus_eh_e$.
Notice that for any $w\in \calw(\G)$ the natural projection of $T_w$ to $E_w$
induces a $m_w$--covering of $F_w\to E_w$ with Galois group $\bfz_{m_w}$.
Moreover, the monodromy action $h_w$ can be identified with the action
of the generator $\hat{1}$ of this Galois group. In particular, $h_w$ has
finite order. Now, it is elementary to verify that the other monodromy
operators  $h_a$ ($a\in \cala$) and $h_e$ ($e\in\cale$) are also of finite 
order. On the other hand,  $h^*_{|G|}$ is  finite by its construction.

Let $N_0$ be an integer such 
that $\lambda^{N_0}=1$ for any eigenvalue $\lambda$ of $h_1$ and 
$h_v$.  Then $\delta^0(im(h_1^{N_0}-1))= 0$, hence
$H^1(|G|)\subset im(h_1^{N_0}-1)$. But, by the above corollary (\ref{9.8})
(part b), 
these spaces have the same dimension. Therefore, we have:

\bekezdes{Corollary.}{9.12}
{\em Fix an $N_0$ such that $\lambda^{N_0}=1$ for any eigenvalue 
$\lambda $ of $h_1(f)$. Then
the pairs $(H^1(|G|),h^*_{|G|})$ and $(im(h_1^{N_0}-1),h_1)$ are 
isomorphic.  This 
identification is compatible with the generalized eigenspace decomposition,
in particular, for any eigenvalue $\lambda$
 there exists the following commutative diagram:

$$\begin{array}{ccccccccc}
0&\to&H^1(|G|)_{\lambda}&\to&H^1(F)_{\lambda}&\to&coker(h_1^{N_0}-1)_{\lambda}&\to&0\\
&&\Big\downarrow\vcenter{%
\rlap{$h_{|G|}^*$}}&&\Big\downarrow\vcenter{%
\rlap{$h_1$}}&&\Big\downarrow\vcenter{%
\rlap{$h_1$}}&&\\
0&\to&H^1(|G|)_{\lambda}&\to&H^1(F)_{\lambda}&\to&coker(h_1^{N_0}-1)_{\lambda}&\to&0.
\end{array}$$
 This shows that for any $\lambda$:
$$\#^2_{\lambda}(f)=H^1(|G|)_{\lambda}.$$
}

\noi
Actually, in this diagram  we can also see all the cohomology groups
$H^1(mod_N(G)|)$ for arbitrary integer $N$. Indeed,
the ``orbit projection'' $o_N:|G|\to |mod_N(G)|$
induces an injective morphism $o_N^*:H^1(|mod_N(G)|)\to H^1(|G|)$
(where $im(o_N^*)$ are the invariant cocycles), therefore $o_N^*$ can be 
identified with the inclusion $$H^1(|mod_N(G)|)=\ker(h_{|G|}^N-1)
\hookrightarrow H^1(|G|).$$

\bekezdes{}{9.13} 
The fact that the $2\times 2$--Jordan blocks can be characterized by a 
graph with a $\bfz$--action, has the following interesting consequence
(cf. \cite{NSM}):

\bekezdes{Proposition.}{9.14}
{\em Assume that $c_{\G}=0$ and $G(X,f)$ is connected. Assume that $h_1$ has
 a Jordan block of size 2 with an eigenvalue $\beta$ with
$\beta^k=1$. Write $k=pq$ with
$(p,q)=1$. Then $h_1$ has a Jordan block of size 2 with some 
eigenvalue, say  $\lambda$, 
such that either $\lambda^p=1$ or $\lambda^q=1$.}

\vspace{2mm}

\noindent 
{\em Proof.}\ Fix a vertex $w^*\in\calw(\G)$ such that
there is an arrow $a\in\cala(\G)$ with $(a,w^*)\in\cale(\G)$. Let 
$\cale_n$ be the set of edges of $\G$ connecting two non--arrowhead
vertices. For any $e\in\cale_n$ let $w(e)$ be the vertex of $e$
which has  larger distance in $\G$ from $w^*$ (i.e. $w(e)$ and $w^*$ are
in two different components of $\G\setminus \{e\}$). Then $e\mapsto w(e)$
defines a bijection $\cale_w\to \calw\setminus\{w^*\}$. By
Euler--characteristic argument:
$$1-c_{mod_k(G)}=\sum_{w\in\calw}g.c.d.(k,n_w)-\sum_{e\in\cale_n}
g.c.d.(k,n_e).$$
But $n_{w^*}=1$ (cf. \ref{5.7}), hence:
$$c_{mod_k(G)}=\sum_{e\in\cale_n}g.c.d.(k,n_e)-g.c.d.(k,n_{w(e)}),$$
where $n_{w(e)}|n_e$.

If $\#^2_{\beta}(f)\not=0$ for some $\beta$ with $\beta^k=1$, then 
$c_{mod_k(G)}\not=0$, hence there exists $e\in\cale_n$ with
$(k,n_e)>(k,n_{w(e)})$. But then, for the same $e$, a similar strict 
inequality holds  either
for $p$ or for $q$ (instead of $k$). This ends the proof.\hfill $\diamondsuit$

\bekezdes{Remark.}{9.15} 
The above proof contains a generalization of the result
(\ref{9.7}) of Neumann.
Namely, even if $g\not=0$ (but with $c_{\G(X)}=0$), there exist positive integers
$n_i$ and $m_i$ ($i\in I$), with $m_i|n_i$, such that:
$$\sum_{\lambda^N=1}\#^2_{\lambda}(f)=\sum_ig.c.d.(N,n_i)-g.c.d.(N,m_i).$$

\bekezdes{The finiteness of $h_1$ revisited.}{9.16}
Above we proved the characterization
 $\sum_{\lambda}\#^2_{\lambda}(f)=c_{G(X,f)}$, but we  can  still
ask: when is this number zero?

For example, L\^{e} D. T. proved  that the monodromy of an irreducible
germ $f:(\bfc^2,0)\to (\bfc,0)$ is always finite \cite{Le}.

In L\^{e}'s result the irreducibility assumption is really important.
Indeed, if we take
$f:(\bfc^2,0)\to (\bfc,0)$  given by $f(x,y)=(x^2+y^3)(x^3+y^2)$, then by
(\ref{5.12}), the resolution of the double covering $\bfc^2_{f,2}$ has a cycle,
hence $h_1(f)$ is not finite (historically, this is the first germ with
non--finite monodromy, found by A'Campo \cite{A'CMon}). 

On the other hand,  if $f:(X,x)\to (\bfc,0)$ has a 
finite monodromy $h_1(f)$, then 
$c_{\G(X)}$ must be zero (cf. \ref{9.5}).

\vspace{2mm}

So, this motivates to investigate 
the case when $c_{\G(X)}=0$ and $\#\cala(\G)=1$. Under a
more restrictive condition, when $L_X$ is an {\em integer} homology sphere,
W. Neumann (using a stronger version of the identity (\ref{9.7})) proved
that  $h_1(f)$ is finite, provided that $\#\cala(\G)=1$ (see \cite{EN}
page 111, or \cite{Ne2}). The next result is a generalization of the 
results of L\^{e} and Neumann.

\bekezdes{Theorem.}{9.17} (cf. also with \cite{NSM})
{\em Assume that $L_X$ is a rational homology sphere
and $f:(X,x)\to (\bfc,0)$ defines an isolated singularity with $\#\cala(\G)=1$.
Then if $\#^2_{\lambda}\not=0$ for some $\lambda$ of order $k$, then
there exists a prime number $p\geq 2$ such that $p|k$ and $p^2|\ |H_1(L_X,\bfz)|$.

In particular,  if $|H_1(L_X,\bfz)|$ is square free, then 
$h_1(f)$ is finite for any $f$.

(Recall that the order $|H_1(L_X,\bfz)|$ of the torsion group $H_1(L_X,\bfz)$
is exactly $|\det (E_w\cdot E_v)|$.)}

\vspace{2mm}

\noindent {\em Proof.}\  
Assume that $\#^2_{\lambda}\not=0$ with $\lambda^k=1$.
Then by (\ref{9.14}), there is  a prime number $p$ and eigenvalue $\eta$
such that $p|k$, $\eta^{p^{\alpha}}=1$ for some $\alpha\geq 1$, and 
$\#^2_{\eta}\not=0$. Write  the characteristic polynomial  $\Delta(t)=
\det (th_1-1)$ as a product of cyclotomic  polynomials 
$\prod_l\phi_l(t)^{t_l}$.
By Wang's exact sequence and Alexander duality, $\coker
(h_1-1)=H_1(L_X,\bfz)$. 
Then $\Delta(1)=\prod_l\phi_l(1)^{t_l}=|H_1(L_X,\bfz)|$. Now, if $l$ is not 
a prime power, then $\phi_l(1)=1$, whereas $\phi_l(1)=q$ if $l$ is a power of
 the prime $q$. Since $h_1$ has a double eigenvalue of order a power of 
the prime $p$, we obtain the divisibility $p^2|\Delta(1)$.\hfill$\diamondsuit$

\bekezdes{Example.}{ns1}\cite{NSM}\ Consider the embedded resolution  graph of
$f:(X,x)\to (\bfc,0)$, where $(X,x)=(\{x^{10}+x^3y^2+y^3+z^2=0\},0)\subset
(\bfc^3,0)$ and $f(x,y,z)=x$:

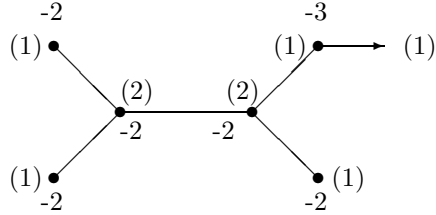
\begin{picture}(200,100)(0,0)
\put(50,25){\circle*{4}}
\put(50,75){\circle*{4}}
\put(75,50){\circle*{4}}
\put(125,50){\circle*{4}}
\put(150,25){\circle*{4}}
\put(150,75){\circle*{4}}
\put(50,75){\line(1,-1){25}}
\put(50,25){\line(1,1){25}}
\put(75,50){\line(1,0){50}}
\put(125,50){\line(1,1){25}}
\put(125,50){\line(1,-1){25}}
\put(150,75){\vector(1,0){25}}
\put(33,72){(1)}
\put(33,22){(1)}
\put(75,55){(2)}
\put(115,55){(2)}
\put(133,72){(1)}
\put(155,22){(1)}
\put(182,72){(1)}
\put(45,85){-2}
\put(45,13){-2}
\put(75,40){-2}
\put(110,40){-2}
\put(145,13){-2}
\put(145,85){-3}
\end{picture}

Then $|\det(E_v\cdot E_w)|=|H_1(L_X,\bfz)|=4$. 
Then $G$ can be computed from \ref{5.7}, hence
the monodromy $h_1$ contains exactly
one Jordan block of size 2 and the eigenvalue of this block is $=-1$. 

\vspace{2mm}

Finally, we  ask:  in the above result (\ref{9.17}), 
is it really important for $L_X$ to be 
a rational homology sphere? Is the monodromy finite, for example, 
if $c_{\G(X)}=0$, \ $\#\cala(\G(X,f))=1$ and Tors$(H_1(L_X,\bfz))=0$ (i.e. the 
intersection matrix $(E_w\cdot E_v)$ is unimodular)? The answer is negative!
A possible counterexample is the following:

\bekezdes{Example.}{9.18} In theorem (\ref{9.17}) not 
only $c_{\G}=0$ is important but 
also $g=0$. To see this, take 
 $(X,x)=(\{x^2+y^7-z^{14}=0\},0)\subset
(\bfc^3,0)$ and the function 
$f_2(x,y,z)=z^2-y$ as in (\ref{3.17}) and (\ref{5.10}).
Let $P$ be the intersection point of the strict transform $S_a$ of
$\{f_2=0\}$ with the exceptional divisor $E$. Then, in some local  coordinate
system $(u,v)$ of $P$, $\{u=0\}$ represents $E$ (in a neighborhood of $P$),
$\{v=0\}$ represents $S_a$, and $f_2=u^2v^2$
(cf. \ref{3.17}). Consider $g=z$. Since $z$ in the
neighborhood of $P$ can be represented as $z=u$ (modulo a local invertible
germ), $f_2+g^k$ near $P$ has the form $u^2v^2+u^{k}$. For example, if $k=3$,
then one needs two more blowing ups in order to resolve $f_2+g^k$.
In this case, the graph $\G(X,-y+z^2+z^3)$ and its covering $G(X,-y+z^2+z^3)$
are the following:

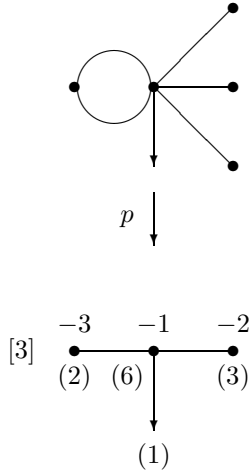
\begin{picture}(400,200)(100,-40)
\put(270,20){\circle*{4}}
\put(210,20){\circle*{4}}
\put(240,20){\circle*{4}}
\put(210,20){\line(1,0){60}}

\put(240,20){\vector(0,-1){30}}
\put(190,20){\makebox(0,0){$[3]$}}
\put(210,10){\makebox(0,0){$(2)$}}
\put(210,30){\makebox(0,0){$-3$}}
\put(230,10){\makebox(0,0){$(6)$}}
\put(240,30){\makebox(0,0){$-1$}}
\put(240,-20){\makebox(0,0){$(1)$}}

\put(270,10){\makebox(0,0){$(3)$}}
\put(270,30){\makebox(0,0){$-2$}}

\put(240,80){\vector(0,-1){20}}
\put(230,70){\makebox(0,0){$p$}}

\put(270,150){\circle*{4}}
\put(270,120){\circle*{4}}
\put(270,90){\circle*{4}}

\put(240,120){\line(1,1){30}}
\put(240,120){\line(1,0){30}}
\put(240,120){\line(1,-1){30}}

\put(210,120){\circle*{4}}
\put(240,120){\circle*{4}}
\put(225,120){\circle{29}}
\put(240,120){\vector(0,-1){30}}

\end{picture}
 
So, in this case, the intersection matrix of
$\G$ is unimodular, hence $H_1(L_X,\bfz)$
is torsion free, $c_{\G}=0$, and $\#{\cala}=1$, but $\#^2_{-1}=1$, hence
$h_1$ is not finite. The reason is that $g\not=0$. 

This also shows  that  part b) of Theorem B in \cite{NSM} is true only
with the additional assumption $g=0$. (The other statements of 
\cite{NSM} are correct; the error in (B,b) comes from the identity (1) on 
page 591.)

{

}

\end{document}